\documentclass[bj,authoryear,noshowframe]{imsart}

\usepackage{multirow}
\usepackage{paralist}
\usepackage{bm,bbm}
\usepackage{stackengine}
\stackMath
\newcommand\tenq[2][1]{%
\def\useanchorwidth{T}%
\ifnum#1>1%
\stackunder[0pt]{\tenq[\numexpr#1-1\relax]{#2}}{\scriptscriptstyle\thicksim}%
\else%
\stackunder[1pt]{#2}{\scriptscriptstyle\thicksim}%
\fi%
}
\usepackage{mathtools}

\startlocaldefs
\theoremstyle{plain}

\newtheorem{theorem}{Theorem}[section]
\newtheorem{lemma}[theorem]{Lemma}
\newtheorem{proposition}[theorem]{Proposition}
\theoremstyle{remark}
\newtheorem{assumption}[theorem]{Assumption}

\newtheorem{remark}[theorem]{Remark}

\endlocaldefs

\newcommand{\cqfd}{\ensuremath{\hfill\Box}}

\newcommand{\pr}{^{\prime}}
\newcommand{\n}{^{(n)}}
\newcommand{\m}{\!\!\!\!\!\!}
\newcommand{\npr}{^{(n)\prime}}
\newcommand{\rP}{{\rm P}}
\newcommand{\rd}{{\rm d}}
\newcommand{\Xb}{{\mathbf{X}}}
\newcommand{\tZ}{\mathbf{Z}}
\newcommand{\pms}{{\scriptscriptstyle \pm}}

\allowdisplaybreaks

\usepackage{pdfpages}

\begin{document}

\begin{frontmatter}
\title{{Distribution-free} tests of multivariate independence based on center-outward quadrant, Spearman, Kendall, and {van~der~Waerden} statistics}
\runtitle{{Distribution-free} tests of multivariate independence}

\begin{aug}
\author[A]{\inits{H.}\fnms{Hongjian}~\snm{Shi}\ead[label=e1]{hongjian.shi@tum.de}\orcid{0000-0001-5552-2082}}
\author[A]{\inits{M.}\fnms{Mathias}~\snm{Drton}\ead[label=e2]{mathias.drton@tum.de}\orcid{0000-0001-5614-3025}}
\author[B]{\inits{M.}\fnms{Marc}~\snm{Hallin}\ead[label=e3]{mhallin@ulb.ac.be}\orcid{0000-0002-6599-7409}}
\author[C]{\inits{F.}\fnms{Fang}~\snm{Han}\ead[label=e4]{fanghan@uw.edu}\orcid{0000-0003-2996-5693}}
\address[A]{Department of Mathematics, TUM School of Computation, Information and Technology, Technical University of Munich, 85748 Garching bei M\"unchen, Germany,
\printead{e1,e2}}

\address[B]{ECARES and Department of Mathematics, Universit\'e Libre de Bruxelles, Brussels, Belgium,
\printead{e3}}

\address[C]{Department of Statistics, University of Washington, Seattle, WA 98195, USA,
\printead{e4}}
\end{aug}

\begin{abstract}
Due to the lack of a canonical ordering in ${\mathbb R}^d$ for $d>1$, defining multivariate generalizations of the classical univariate ranks has been a long-standing open problem in statistics.  Optimal transport has been shown to offer a solution in which multivariate ranks are obtained by transporting data points to a grid that approximates a uniform reference measure \citep{MR3611491, hallin2017distribution, MR4255122}, thereby inducing ranks, signs, and a data-driven ordering of ${\mathbb R}^d$.  We take up this new perspective to define and study multivariate analogues of the sign covariance/quadrant statistic, Spearman's rho, Kendall's tau, and van der Waerden covariances.  The resulting tests of multivariate independence are fully distribution-free, hence uniformly valid irrespective of the actual (absolutely continuous) distribution of the observations.  Our results provide the asymptotic distribution theory for these new test statistics, with asymptotic approximations to critical values to be used for testing independence between random vectors, as well as a power analysis of the resulting tests in an extension of the so-called (bivariate) Konijn model.  This power analysis includes a multivariate Chernoff--Savage property gua\-ranteeing~that, under elliptical generalized Konijn models, the asymptotic relative efficiency of our van der Waerden tests with respect to Wilks' classical (pseudo-)Gaussian procedure is strictly larger than or equal to one, where equality is achieved under Gaussian distributions only.  We similarly provide a lower bound for the asymptotic relative efficiency of our Spearman procedure with respect to Wilks\!' test, thus extending the classical result by Hodges and Lehmann on the asymptotic relative efficiency, in univariate location models, of Wilcoxon tests with respect to the Student ones.
\end{abstract}


\end{frontmatter}

\section{Introduction}

Testing independence between two observed random variables $X_1$ and $X_2$ is a fundamental problem in statistical inference and has important applications, basically, in all domains. In a bivariate context, when both $X_1$ and $X_2$ are univariate Gaussian, the classical tests are based on empirical correlations.  These tests remain asymptotically valid under non-Gaussian distributions with finite variances and, therefore, are widely used in practice as pseudo-Gaussian tests. Rank-based alternatives, that do not require any moment assumptions, have been proposed at a very early stage: the Spearman and  Kendall correlation coefficients (\citealp{Spearman1904}, \citealp{Kendall1938}), actually, were among the first applications of rank-based methods in statistical inference, well before \citet{10.2307/3001968} gave his rank sum and signed rank tests for location. The general opinion, however, was that the (high) price to be paid for the extended validity of these rank tests was their poor performance in terms of power ... until \citet{MR79383} and \citet{MR100322} disproved this fact by showing that, in univariate location models, the asymptotic relative efficiency (ARE) of Wilcoxon tests relative to Student's $t$ ones never falls below 0.864 while the same ARE, for the van der Waerden (normal score) version of rank tests is uniformly not smaller than one. Originally proved in the context of univariate location models, these Hodges--Lehmann and Chernoff--Savage results have been extended to the context of bivariate independence and the so-called Konijn model  (\cite{MR79384}; see also Chapter~III.6.1 in \cite{MR0229351}) by~\citet{MR2462206}.
 
In this paper, we consider the problem of testing {\it multivariate} independence, that is, independence between random vectors ${\bf X}_1$ and ${\bf X}_2$ with dimensions $d_1$ and $d_2$ such that max$(d_1,d_2)>1$. The classical pseudo-Gaussian test for this problem is Wilks' Gaussian likelihood ratio test \citep{Wilks1935}. Based on the concepts of multivariate {\it center-outward ranks} and {\it signs} recently introduced by \citet{MR4255122}, we construct fully distribution-free tests of the quadrant, Spearman, Kendall, and van der Waerden type and assess their asymptotic performance against generalized versions of the Konijn alternatives.  We also provide (against elliptical generalized Konijn alternatives)  a   Hodges--Lehmann lower bound for the ARE of our Spearman tests, and establish a Chernoff--Savage property for our van der Waerden tests, which makes Wilks' pseudo-Gaussian ones Pitman-nonadmissible in this context.

\subsection{Testing multivariate independence}\label{subsec:tmi}

The problem of testing independence between two random vectors ${\bf X}_1$ and ${\bf X}_2$ with dimensions $d_1$ and~$d_2$ and unspecified densities is significantly harder for max$(d_1,d_2)>1$ than for $d_1=d_2=1$---due, mainly, to the difficulty of defining a multivariate counterpart to univariate ranks. The first attempt to provide a rank-based alternative to the Gaussian likelihood ratio method of \citet{Wilks1935} was developed in Chapter~8 of \citet{MR0298844} and, for almost thirty years, has remained the only rank-based approach to the problem. The proposed tests, however, are based on componentwise ranks, which are not distribution-free---unless, of course, both vectors have dimension one, in which case we are back to the traditional context of bivariate independence (e.g., Chapter~III.6 of \citeauthor{MR0229351}, \citeyear{MR0229351}). This issue persists in more recent work, e.g., that of \citet{lin2017copula}, \citet{MR3842884}, and \citet{MR4374647}. Other types of multivariate independence tests include \citet{MR0067429}, \citet{MR696054}, \citet{MR2255909}, \citet{MR2382665}, \citet{MR2956796}, \citet{MR3068450}, \citet{MR4352523}, and \citet{MR4589063}, to name a few; see \citet{MR4399094} and \citet{MR4474478} for a more complete review.

Alternatives to the \citeauthor{MR0298844} test have appeared with the developments of multivariate concepts of signs and ranks. Based on \citet{MR1134492}'s concept of  {\it interdirections}, \citet{MR2691505} and \citet{MR1467849} proposed a sign test extending the univariate quadrant test \citep{MR39190}. \citet{MR2201019} also proposed, based on the so-called standardized {\it spatial signs}, a sign test which, under elliptical symmetry assumptions, is asymptotically equivalent to the Gieser and Randles one. {\it Spatial ranks} are introduced, along with {\it spatial signs}, in \citet{MR2201019}, where multivariate analogs of Spearman's {\it rho} and Kendall's {\it tau}  are considered; the Spearman tests (involving Wilcoxon scores) are extended, in \citet{MR2088309}, to the case of arbitrary square-integrable score functions, which includes van der Waerden (normal score) tests. All these multivariate rank-based tests are enjoying, under elliptical symmetry, many of the attractive properties of their traditional univariate counterparts. In particular, under the assumptions of elliptical symmetry, they are asymptotically distribution-free. Local powers are obtained in \citet{MR2088309} against elliptical extensions of the  {\it Konijn alternatives} \citep{MR79384}.  Hodges--Lehmann and Chernoff--Savage results have been established by \citet{MR2462206}; the latter entail the Pitman-nonadmissibility, in this elliptical context, of Wilks' classical procedure.  

The validity of all these tests, however, is limited to subclasses of distributions---essentially, the fa\-mi\-ly of elliptically symmetric ones---which, from the perspective we are taking here, are overly restrictive; the assumption of elliptical symmetry, in particular, is extremely strong, and unlikely to hold in most applications.  Moreover, there is a crucial difference between finite-sample and asymptotic distribution-freeness.  Indeed, one should be wary that a sequence of tests $\psi\n$ with asymptotic size~$\lim_{n\to\infty}{\rm E}_{\rP}[\psi\n]=\alpha$ under any element~$\rP$ in a class $\mathcal P$ of distributions does not necessarily have asymptotic size~$\alpha$ under unspecified~$\rP\!\in\!\mathcal P$: the convergence of ${\rm E}_{\rP}[\psi\n]$ to $\alpha$, indeed, typically is not uniform over~$\mathcal P$, so that, in general, $\lim_{n\to\infty}\sup_{\rP\in\mathcal P}{\rm E}_{\rP}[\psi\n]\neq \alpha$. This is the case for most pseudo-Gaussian procedures, including Wilks' test. Genuinely distribution-free tests~$\psi\n\!$, where ${\rm E}_{\rP}[\psi\n]$ does not depend on $\rP$, do not suffer that problem, and this is why finite-sample distribution-freeness is a fundamental property.

Building on the concept of center-outward ranks and signs recently proposed by \citet{MR4255122},  \citet{MR4474478} are introducing fully distribution-free rank-based versions of a class of {\it generalized symmetric covariances} which includes, among others, the sophisticated measures of multivariate dependence proposed by \citet{MR2382665, MR3269983}, \citet{MR3842884}, \citet{MR3737307}, and \citet{MR4185814} (none of which is distribution-free). These rank-based generalized symmetric covariances are cumulating the advantages of distribution-freeness with those of ``universal'' consistency. However, they also are relatively complex, with tricky non-Gaussian asymptotics involving (Proposition~5.1 of \citet{MR4474478}) the eigenvalues of an integral equation. Practitioners, therefore, may prefer simpler, more familiar and easily interpretable extensions of the classical bivariate quadrant, Spearman, Kendall, and van der Waerden tests.  Defining such extensions and studying their performance is the objective of this paper.

\subsection{Center-outward signs and ranks}

For dimension $d>1$, the real space $\mathbbm{R}^d$ lacks a canonical ordering. 
As a result, defining, in dimen\-sion~$d>1$, concepts of signs and ranks enjoying the properties that make the traditional ranks successful tools in univariate inference has been an open problem for more than half a century. One of the most important properties is exact distribution-freeness (for i.i.d.~samples from absolutely continuous distributions).  In an important new development involving optimal transport, the concept of center-outward ranks and signs was proposed recently by \citet{MR3611491}, \citet{hallin2017distribution}, and \citet{MR4255122} and, contrary to earlier concepts such as marginal ranks \citep{MR0298844}, spatial ranks~\citep{MR2598854}, depth-based ranks \citep{MR1212489, MR2329471}, or Mahalanobis ranks and signs \citep{MR1963662, MR1926170}, enjoys a property of {\it maximal ancillarity}---intuitively, {\it ``maximal distribution-freeness''}---under i.i.d.~observations with unspecified absolutely continuous distribution. 

Center-outward ranks and signs are based on measure transportation to the unit ball equipped with a spherical uniform distribution.  They have been used, quite successfully, in a variety of statistical problems: rank tests and R-estimation for VARMA models  \citep{MR4436323, MR4497246},  rank tests for multiple-output regression and MANOVA \citep{MR4646617}, and nonparametric multiple-output quantile regression \citep{del2022nonparametric}.   
  
Other authors (\citealp{MR3485957, MR3698117}, \citealp{MR4404927}, \citealp{MR4571116}, and \citealp{deb2021efficiency}) are considering measure transportation to the Lebesgue uniform over the unit~cube rather than the spherical uniform over the unit ball;  the finite-sample impact of that choice has been studied in \citet{hallin2023finitesample} who show that it is relatively modest. The inverse of such a transport (a quantile function), however, does not preserve the symmetry features of the underlying distribution (for instance, the quantile function of a spherically symmetric distribution is not spherically symmetric), while the quantile function resulting from a transport to the unit ball does.

\subsection{A motivating example}\label{subsec:motex}

The importance and advantages of center-outward sign- and rank-based tests of independence are\linebreak illustrated with the following stock market data example.  The dataset, collected from Yahoo!~Finance (\url{finance.yahoo.com}), contains prices for stocks in the Standard \& Poor's 100 (S\&P 100) index for the years 2003--2012. We will analyze the daily log returns of adjusted closing prices based on~$n=40$ observations---one from each quarter. The time lag between two consecutive observations is three months, in a financial series where~auto\-correlations typically are small (see Figure~S.1 in the supplement). Treating these observations as independent and identically distributed, thus, is unlikely to affect the validity of the tests. 

\begin{figure}[!htbp]
\centering
\includegraphics[width=\textwidth, trim={.35in .35in .35in .35in},clip]{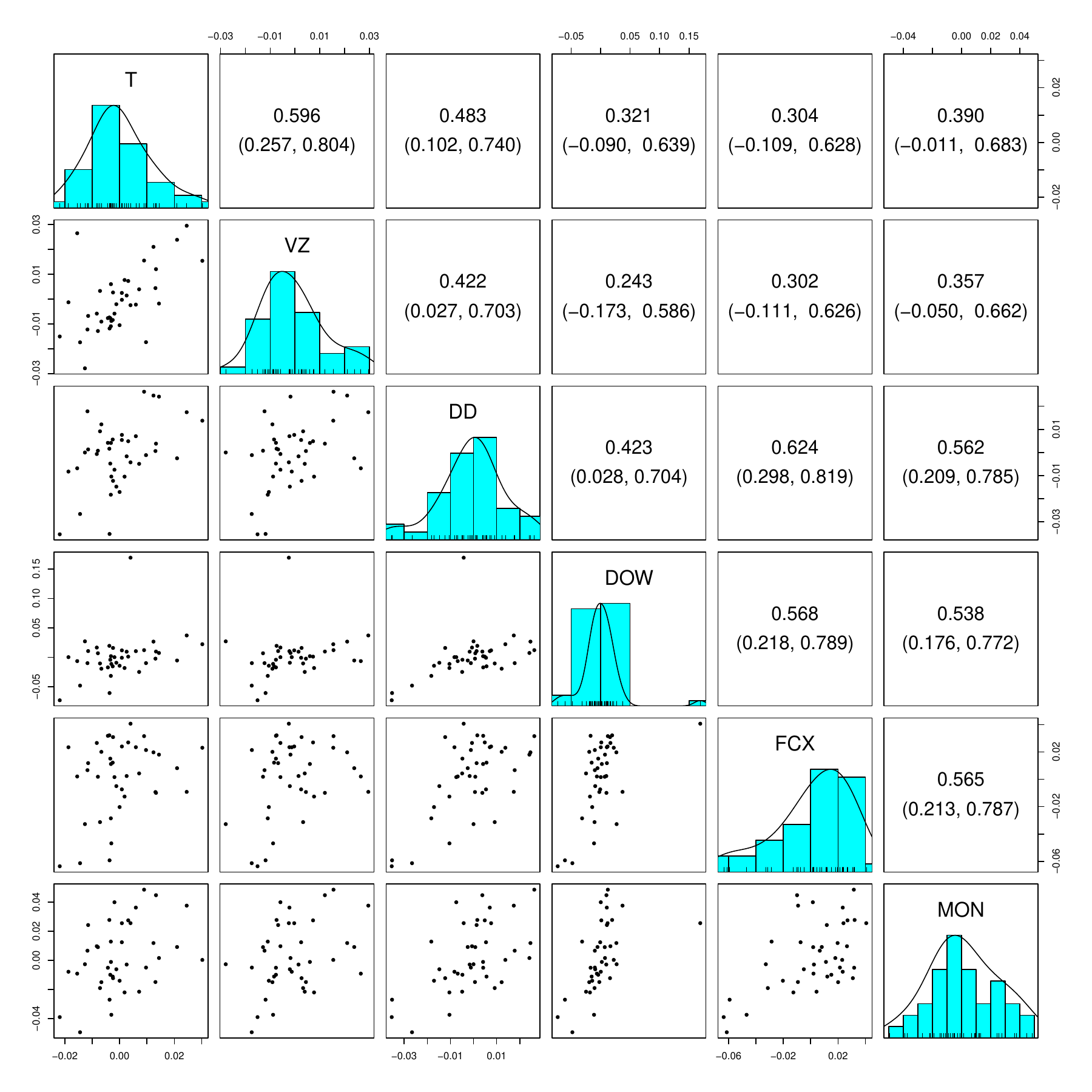}\vspace{-1em}
\caption{Univariate histograms (on the diagonal), bivariate scatter plots (below the diagonal), and the {lag-zero} Pearson  correlation coefficients {with 99\% confidence intervals} (above the diagonal) for the forty observations of daily log returns of six stocks:  ``AT\&T Inc [T]'', ``Verizon Communications [VZ]'', ``Du Pont (E.I.) [DD]'', ``Dow Chemical [DOW]'',  ``Freeport-McMoran Cp \& Gld [FCX]'',  and ``Monsanto Co. [MON]''.\vspace{-1em}}\label{fig:rplot}
\end{figure}

The S\&P 100 stocks are classified into ten sectors by Global Industry Classification Standard (GICS). To illustrate the advantages of the proposed tests, we focus on detecting dependence between two such sectors:  
(1) Telecommunications, including ``AT\&T Inc [T]'' and ``Verizon Communications [VZ]'' stocks;  and (2) Materials, including  ``Du Pont (E.I.) [DD]'', ``Dow Chemical [DOW]'', ``Freeport-McMoran Cp \& Gld [FCX]'', and ``Monsanto Co.\ [MON]'' stocks.   

Figure~\ref{fig:rplot} presents univariate histograms and bivariate scatterplots for the {forty observations (separated by three months) of the daily log returns} of the six stocks mentioned above.  We notice that all univariate marginals are skewed and/or heavy-tailed.  We also notice that the lag-zero Pearson correlations between any stock in \{T,~VZ\} and any stock in \{DD,~DOW,~FCX,~MON\} are relatively small, and that most correlations are non-significant at significance level $1\%$.

We applied the four rank-based tests (quadrant (sign), Spearman, Kendall, van der Waerden) proposed in Section~\ref{testprocsec} to the forty bivariate observations of~(T,~VZ) and either (DD,~DOW), (DD,~FCX),  (DD,~MON),  (DOW,~FCX),  (DOW,~MON), or (FCX,~MON).  These tests are based on the matrices defined in~\eqref{tildeW:sign}--~\!\eqref{scoreW}.  As a benchmark, we also include Wilks' Gaussian likelihood ratio test (LRT).  The~$p$-values for all these tests are reported in Table~\ref{tab:pvalue-real}. One observes that the rank-based tests perform as well as Wilks' LRT in the first three columns while, in the last three, they yield strong evidence against independence whereas Wilks is unable to reject independence at 5\% significance level.  Rank-based sign, Spearman, Kendall, and van der Waerden tests, thus, are able to detect dependencies that Wilks' traditional pseudo-Gaussian test cannot.

{
\renewcommand{\tabcolsep}{2pt}
\renewcommand{\arraystretch}{1.00}
\begin{table}[!htb]
\caption{P-values of the quadrant, Spearman, Kendall, van der Waerden, and Wilks tests of the hypotheses of independence between the bivariate vector  (T,VZ)  and each of the six bivariate vectors (DD,~DOW), 
(DD,~FCX),  (DD,~MON),  (DOW,~FCX),  (DOW,~MON), or (FCX,~MON).\smallskip}
\label{tab:pvalue-real}
\centering
\begin{tabular}{cccccccc}
\hline
                      &    & (DD, DOW)  & (DD, FCX)  & (DD, MON)  & (DOW, FCX) & (DOW, MON) & (FCX, MON) \\
\hline
  quadrant (sign)     &    & $   0.002$ & $   0.050$ & $   0.027$ & $   0.013$ & $   0.008$ & $   0.023$ \\
  Spearman            &    & $\m<0.001$ & $   0.005$ & $   0.004$ & $   0.007$ & $   0.001$ & $   0.011$ \\
  Kendall & \ \ (T,VZ)\ \  & $   0.002$ & $   0.013$ & $   0.016$ & $   0.011$ & $   0.002$ & $   0.036$ \\
  van der Waerden     &    & $\m<0.001$ & $   0.012$ & $   0.011$ & $   0.006$ & $\m<0.001$ & $   0.010$ \\
  Wilks (LRT)         &    & $   0.011$ & $   0.016$ & $   0.009$ & $   0.174$ & $   0.072$ & $   0.074$ \\
\hline
\end{tabular}
\end{table}
}

\subsection{Outline of the paper}

The paper is organized as follows. Section~\ref{sec:2} briefly reviews the notions of center-outward ranks and signs.  Section~\ref{sec:3} introduces our rank tests of multivariate independence.  In Section~\ref{sec:4}, we establish a Chernoff--Savage property and a Hodges--Lehmann result for our center-outward van der Waerden and Spearman tests, respectively.  We briefly conclude in Section~\ref{sec:6}.  Due to the importance, all the proofs are given in Section~\ref{sec:7}.  The auxiliary results and numerical experiments are deferred to the supplement \citep{shi2023distribution-supp}.

\section{Center-outward distribution functions, ranks, and signs} \label{sec:2}

\subsection{Definitions}\label{FQsec}

Denoting  by~${\mathbb{S}_d}$ and ${\mathcal{S}_{d-1}}$, respectively,  the open unit ball and the unit hypersphere in ${\mathbb R}^d$, let ${\rm U}_d$ stand for the spherical \footnote{Namely, ${\rm U}_d$ is the spherical distribution with uniform (over $[0,1]$) radial density---equivalently, the product of a uniform over the distances to the origin and a uniform over the unit sphere ${\cal S}_{d-1}$. For~$d=1$, ${\rm U}_1$ coincides with the Lebesgue uniform over~$(-1,1)$.} uniform distribution over~${\mathbb{S}_d}$.  Let ${\cal P}_d$ denote the class of Lebesgue-absolutely continuous distributions over $\mathbbm{R}^d$.  For any $\rP$ in ${\cal P}_d$, the main result in \citet{MR1369395} then implies the existence of a $\rP$-a.s.\ unique gradient $\nabla\phi$ of a convex (and lower semi-continuous) function~$\phi:\mathbbm{R}^d\to\mathbbm{R}$ such that $\nabla\phi$ pushes $\mathrm{P}$ forward to~${\rm U}_d$, i.e., $\nabla\phi({\bf Z})\sim {\rm U}_d$ under ${\bf Z}\sim\mathrm{P}$.  Call {\it center-outward distribution function} of $\rm P$ any version ${\bf F}_\pms$ of this a.e.~unique gradient. 

Turning to sample versions, denote by $\tZ\n\!\coloneqq \big(\tZ_1\n,\dots, \tZ_n\n\big)$, $n\in~\!\mathbb{N}$ a triangular array of i.i.d.\ $d$-dimensional random vectors with distribution~$\mathrm{P}$.  The {\it empirical center-outward distribution function}~${\bf F}_\pms\n$ of~$\tZ\n$ maps the $n$-tuple $\tZ_1\n,\dots, \tZ_n\n$ to a ``regular'' grid $\mathfrak{G}_n$ of the unit ball~${\mathbb S}_d$. 

This regular grid $\mathfrak{G}_n$ is expected to approximate the spherical uniform distribution ${\rm U}_d$ over the unit ball~${\mathbb{S}_d}$. 
The only mathematical requirement needed for the asymptotic results below is the weak convergence, as~$n\to\infty$, of the uniform discrete distribution over~$\mathfrak{G}_n$ to the spherical uniform distribution~${\rm U}_d$. A spherical uniform i.i.d.~sample of $n$ points over~${\mathbb{S}_d}$ (almost surely) satisfies such a requirement. Since the spherical uniform is highly symmetric, further symmetries can be imposed on~$\mathfrak{G}_n$, though, which only can improve the convergence to ${\rm U}_d$.  In the sequel,  we throughout assume that the grids $\mathfrak{G}_n$ are symmetric with respect to the origin, i.e.,  ${\bf u}\in \mathfrak{G}_n$ implies~${-\bf u}\in\mathfrak{G}_n$, which considerably simplifies formulas. More symmetry, however, can be assumed: see   Section~\ref{Remgrid} below. It should be insisted, however, that imposing such symmetry assumptions does not restrict the generality of the results since the construction of the grid is entirely under control. 

The empirical counterpart ${\bf F}_\pms\n$ of ${\bf F}_\pms$ is defined as the bijective mapping from~$\tZ_1\n,\dots, \tZ_n\n$ to the grid~$\mathfrak{G}_n$ that minimizes $\sum_{i=1}^n\big\Vert {\bf F}_\pms\n (\tZ_i\n) - \tZ_i\n \big\Vert ^2$, where $\|\cdot \|$ denotes the Euclidean norm.  That mapping is unique with probability one; in practice, it is obtained via a simple optimal assignment (pairing) algorithm---a linear program; see \citet[Section~5.1]{MR4399094} for a review and references therein.  Call {\it rescaled center-outward rank} of~$\tZ_i\n$ the modulus  
\begin{equation}\label{rankdef}
\widetilde{R}\n_{i;{{\pms}}}\coloneqq \big\Vert {\bf F}_\pms\n (\tZ_i\n)\big\Vert, \quad i=1,\ldots,n
\end{equation} 
and {\it center-outward sign} of~$\tZ_i\n$ the unit vector 
\begin{equation}\label{signdef}
{\bf S}\n_{i;{\pms}}\coloneqq {\bf F}_\pms\n (\tZ_i\n)\big/\big\Vert {\bf F}_\pms\n (\tZ_i\n)\big\Vert\quad \text{for ${\bf F}_\pms\n (\tZ_i\n)\neq{\bf 0}$;}
\end{equation}
put ${\bf S}\n_{i;{{\pms}}}={\bf 0}$ for ${\bf F}_\pms\n (\tZ_i\n) ={\bf 0}$. When the grid $\mathfrak{G}_n$ is obtained, as in Section~\ref{Remgrid} below, with a factori\-zation of $n$ into $n_Rn_S + n_0$, call $R\n_{i;{{\pms}}}\coloneqq (n_R +1)\widetilde{R}\n_{i;{{\pms}}}$, which  takes values\footnote{The value $R\n_{i;{{\pms}}} = 0$ is attained only if $n_0\neq 0.$} $(0), 1,\ldots,n_R$, the {\it center-outward rank} of~$\tZ_i\n$.

\subsection{Grid selection}\label{Remgrid}

The way (rescaled) ranks and signs in \eqref{rankdef} and \eqref{signdef} are constructed, thus, depends on the way the grid~$\mathfrak{G}_n$ is selected.  As already mentioned, the only requirement for asymptotic results (as in Proposition~\ref{H2018} below) is the weak convergence to the spherical uniform ${\rm U}_d$ of the empirical distribution,~${\rm U}\n_{\mathfrak G}$, say, over $\mathfrak{G}_n$. The closer ${\rm U}\n_{\mathfrak G}$ to ${\rm U}_d$, the better. Imposing on $\mathfrak{G}_n$ some of the many symmetries of  ${\rm U}_d$ only can improve the finite-sample performance of ${\rm U}\n_{\mathfrak G}$ as an approximation of ${\rm U}_d$.  Since ${\rm U}_d$ is the product of a uniform over the distances to the origin and a uniform over the unit sphere ${\cal S}_{d-1}$, it is natural to select $\mathfrak{G}_n$ such that ${\rm U}\n_{\mathfrak G}$ similarly factorizes.  This can be obtained as follows: \smallskip 
\begin{compactenum}
\item[{\it (a)}] first factorize $n$ into $n=n_Rn_S + n_0$, with $0\leq n_0<\min(n_R, n_S)$ and $\min(n_R, n_S)\to\infty$ as $n\to\infty$;\footnote{Note that this implies that $n_0/n = o(1)$ as $n\to\infty$. See \citet[Section~7.4]{mordant2021transporting} for a discussion of the selection of  $n_R$ and~$n_S$. }
\item[{\it (b)}] next consider a ``regular array'' $\mathfrak{S}_{n_S}\coloneqq \{{\bf s}^{n_S}_1,\ldots,{\bf s}^{n_S}_{n_S}\}$ of $n_S$ points on the sphere ${\cal S}_{d-1}$ (see Remark~\ref{Remgrid'} below); 
\item[{\it (c)}] construct the grid $\mathfrak{G}_n$ consisting in the collection of the $n_Rn_S  $ gridpoints of the form
$$\big(r/\big(n_R +1\big)\big){\bf s}^{n_S}_s, \quad r=1,\ldots,n_R,~s=1,\ldots,n_S,$$
along with ($n_0$ copies of) the origin in case $n_0\neq 0$: in total $n-(n_0 -1)$ or $n$ distinct points, thus, according as~$n_0>0$ or $n_0=0$.
\end{compactenum}

\begin{remark}\label{Remgrid'}
By ``regular'' array $\mathfrak{S}_{n_S}$ over ${\cal S}_{d-1}$  in  (b) above, we mean ``as regular as possible'' an array~$\mathfrak{S}_{n_S}$---in the sense, for example, of the {\it low-discrepancy sequences} of the type considered in numerical integration, Monte-Carlo methods, and experimental design.\footnote{See also \citet{hallin2023finitesample} for a spherical version of the so-called {\it Halton sequences}.}  The asymptotic results in Proposition~\ref{H2018} below  only require the weak convergence, as~$n_S\to\infty$, of the uniform discrete distribution over~$\mathfrak{S}_{n_S}$  to the uniform distribution over ${\cal S}_{d-1}$. A uniform i.i.d.~sample of points over~${\cal S}_{d-1}$ (almost surely) satisfies that requirement. However, one can easily construct arrays that are ``more regular''  than an i.i.d.~one.  In particular,  in order for the grid  $\mathfrak{G}_n$ described in {\it (a)}--{\it (c)} above to be symmetric with respect to the origin, one could select an even value of $n_S$ and see that~${\bf s}^{n_S}_s\in \mathfrak{S}_{n_S}$  implies $-\, {\bf s}^{n_S}_s\in\mathfrak{S}_{n_S}$, so that~$\sum_{i=1}^n {\bf S}\n_{i;{{\pms}}}=n_R \sum_{s=1}^{n_S}{\bf s}^{n_S}_s =\boldsymbol 0$.
\end{remark}

\begin{remark}\label{Remgrid''}
Some desirable finite-sample properties, such as strict independence between the ranks and the signs, hold with grids constructed as in (a)--(c) above provided, however, that~$n_0=0$ or 1. This is due to the fact that the mapping from the sample to the grid is no longer injective for $n_0\geq~\! 2$. This fact, which  has no asymptotic consequences (since the number~$n_0$ of tied values involved is $o(n)$ as~$n\to\infty$), is easily taken care of by performing the following tie-breaking device in step (c) of the construction of~$\mathfrak{G}_n$: 
{\it (i)} randomly select  $n_0$ directions ${\bf s}^0_1,\ldots,{\bf s}^0_{n_0}$ in~$\mathfrak{S}_{n_S}$, then 
{\it (ii)} replace the $n_0$ copies  of the origin with the new gridpoints 
$[1/2(n_R+1)]{\bf s}^0_1,\ldots,[1/2(n_R+1)]{\bf s}^0_{n_0}.$
This new grid (for simplicity, the same notation ${\mathfrak{G}}_n$ is used as for the original one) no longer has multiple points: the optimal pairing between the sample and the grid is bijective and the resulting (rescaled) ranks and signs are mutually independent.
\end{remark}

\subsection{Main properties}\label{Propsec} 

This section summarizes some of the main properties of the concepts defined in Section~\ref{FQsec}; further properties and the proofs can be found in \citet{MR4255122}, \citet{MR4646617}, and \citet{MR4394914}.  In the following propositions, the grids $\mathfrak{G}_n$ used in the construction of the empirical center-outward distribution function ${\bf F}\n_{\pms}$ are only required to satisfy the minimal assumption of an empirical distribution ${\rm U}\n_{\mathfrak G}$ converging weakly to ${\rm U}_d$.

\begin{proposition}\label{H2018} Let ${\bf F} _{\pms}$ denote the center-outward distribution function of~${\rm P}\in{\cal P}_d$.  Then,
\begin{compactenum}\vspace{.5mm}
\item[(i)]${\bf F} _{\pms}$ is a probability integral transformation of $\mathbbm{R}^d$, that is, ${\bf F} _{\pms}({\bf Z})\sim {\rm U}_d$ if and only if ${\bf Z}\sim {\rm P}$; under~${\bf Z}\sim {\rm P}$, $\Vert{\bf F} _{\pms}({\bf Z})\Vert$ is uniform over $[0, 1)$, ${\bf F} _{\pms}({\bf Z})/\Vert{\bf F} _{\pms}({\bf Z})\Vert$ is uniform over the sphere~${\cal S}_{d-1}$, and they are mutually independent. 
\end{compactenum}\vspace{.5mm}
Let ${\bf Z}\n_1,\ldots ,{\bf Z}\n_n$ be i.i.d.\ with distribution ${\rm P}\in{\mathcal P}_d$ and empirical center-outward distribution function~${\bf F}_\pms\n$. Then, if there are no ties in the grid $\mathfrak{G}_n$, 
\begin{compactenum}\vspace{.5mm}
\item[(ii)] $\big({\bf F}\n_{\pms}({\bf Z}\n_{1}),\ldots , {\bf F}\n_{\pms}({\bf Z}\n_{n}) \big)$ is uniformly distributed over the $n!$ permutations of~$\mathfrak{G}_n$;
\item[(iii)] the $n$-tuple $\big({\bf F}\n_{\pms}({\bf Z}\n_{1}),\ldots , {\bf F}\n_{\pms}({\bf Z}\n_{n}) \big)$ is {\em strongly essentially maximal ancillary};\footnote{See Section~2.4 and Appendices D.1 and D.2 of  \citet{MR4255122} for a precise definition  and a proof of this essential property.}
\item[(iv)] (pointwise convergence) for any $i=1,\ldots,n$,
\begin{equation*}
\displaystyle{}\Big\Vert   {\bf F}\n _{\pms}({\bf Z}\n_i) - {\bf F} _{\pms}({\bf Z}\n_i) \Big\Vert  \rightarrow 0 ~\textrm{\em a.s.} \quad \text{as}~n\to\infty.
\end{equation*}
\end{compactenum}\vspace{.5mm}
Assume, moreover, that ${\rm  P}$ is in the so-called class ${\cal P}^+_d\subset{\cal P}_d$ of distributions {\it with nonvanishing densities}---namely, the class of  distributions with density $f\coloneqq {\rm d P}/{\rm d}\mu_d$ ($\mu_d$ the $d$-dimensional Lebesgue measure) such that, for all~$D\in\mathbbm{R}^+$, there exist constants $\lambda^-_{D;\mathrm{P}}$ and  $\lambda^+_{D;\mathrm{P}}$ satisfying 
 \begin{equation}\label{nonvanprop}
 0<\lambda^-_{D;\mathrm{P}}\leq f({\bf z})
\leq \lambda^+_{D;\mathrm{P}}<\infty
\end{equation}
for all $\bf z$ with $\Vert{\bf z}\Vert \leq D$. 
Then, 
\begin{compactenum}\vspace{.5mm}
\item[(v)] there exists 
a version of ${\bf F}_\pms$  defining a homeo\-morphism between the punctured unit ball~${\mathbb S}_d\!\setminus~\!\{{\bf 0}\}$ and $\mathbbm{R}^d\setminus {\bf F}_\pms^{-1}(\{{\bf 0}\})$; that version  has a continuous inverse ${\bf Q}_\pms$ (with domain ${\mathbb S}_d\!\setminus\!\{{\bf 0}\}$), which naturally  qualifies as~$\rP$'s {\em center-outward quantile function}; 
\item[(vi)] (Glivenko--Cantelli) 
\begin{equation*}
\displaystyle{\max_{1\leq i\leq n}}\Big\Vert   {\bf F}\n _{\pms}({\bf Z}\n_i) - {\bf F} _{\pms}({\bf Z}\n_i) \Big\Vert  \longrightarrow 0 \ \textrm{\em a.s.} \quad \text{as}~n\to\infty.
\end{equation*}
\end{compactenum}
\end{proposition}

Results {\it (v)} and {\it (vi)} are due to \citet{MR3886582}, \citet{hallin2017distribution}, \citet{del2018smooth}, \citet{MR4255122} and  can be extended \citep{MR4147635}
to a more general\footnote{Namely, ${\cal P}_d^+\subsetneq{\cal P}_d^{\rm conv}\subsetneq{\cal P}_d$.} class~${\cal P}_d^{\rm conv}$ of absolutely continuous distributions, with density $f$ supported on a open convex set of $\mathbbm{R}^d$ but not necessarily the whole space, 
while the definition of ${\bf F}_\pms$ given in \citet{MR4255122} aims at selecting, for each~$\rP\in{\cal P}_d$, a version of  $\nabla\phi$ which, whenever $\rP\in{\cal P}_d^{\rm conv}$, is yielding the homeomorphism mentioned in {\it (v)}. For the sake of simplicity, since we are not interested in quantiles, we stick here to the $\rP$-a.s.\ unique definition given above for $\rP\in{\cal P}_d$. \pagebreak

Center-outward distribution functions, ranks, and signs also inherit, from the invariance of squared Euclidean distances, elementary but quite remarkable invariance and equivariance properties under shifts, orthogonal transformations, and global rescaling (see \citet{MR4646617}). Denote by ${\bf F}^{{\bf Z}}_{\pms}$ the center-outward distribution function of $\bf Z$ and by~${\bf F}^{{\bf Z};\mathfrak{G}_n}_{\pms}$ the empirical center-outward distribution function of an i.i.d.~sample~${\bf Z}_1,\ldots,{\bf Z}_n$ associated with a grid $\mathfrak{G}_n$. 
For $\boldsymbol{\mu}\in\mathbbm{R}^d$, $h\in\mathbbm{R}^+$, and~${\bf O}$ a~$d\times d$ orthogonal matrix, denote
by ${\bf F}^{\boldsymbol{\mu}+ h{\bf O}{\bf Z};{\bf O}\mathfrak{G}_n}_{\pms}$ the empirical center-outward distribution function computed from  the sample~$\boldsymbol{\mu}+ h{\bf O}{\bf Z}_1,\ldots, \boldsymbol{\mu}+ h{\bf O}{\bf Z}_n$ and the grid ${\bf O}\mathfrak{G}_n\coloneqq \{{\bf O}{\bf u}:{\bf u}\in\mathfrak{G}_n\}$. 

\begin{proposition}\label{invF}
Let $\boldsymbol{\mu}\in\mathbbm{R}^d$, $h\in\mathbbm{R}^+$, and ${\bf O}$ be a~$d\times d$ orthogonal matrix. Then,\smallskip
\begin{compactenum}
\item[(i)] ${\bf F}^{\boldsymbol{\mu}+ h{\bf O}{\bf Z}}_{\pms} (\boldsymbol{\mu} + h{\bf O}{\bf z})= {\bf O}{\bf F}^{\bf Z}_{\pms}({\bf z})$, ${\bf z}\in\mathbbm{R}^d$;\smallskip
\item[(ii)] 
$\label{equiv}
{\bf F}^{\boldsymbol{\mu}+ h{\bf O}{\bf Z};{\bf O}\mathfrak{G}_n}_{\pms} (\boldsymbol{\mu} + h{\bf O}{\bf Z}_i)= {\bf O}{\bf F}^{{\bf Z};\mathfrak{G}_n}_{\pms}({\bf Z}_i), ~  i=1,\ldots,n.
$
\end{compactenum}
\end{proposition}

Note that the orthogonal transformations in Proposition~\ref{invF} include the permutations of $\bf Z$'s components. Invariance with respect to such permutations is an essential requirement for hypothesis testing in multivariate analysis.

\section{Rank-based tests for multivariate independence}\label{sec:3}
\subsection{Center-outward test statistics for multivariate independence}
In this section, we describe the test statistics we are proposing for testing independence between two random vectors.  Consider a sample $(\Xb \pr_{11},\Xb \pr_{21})\pr,(\Xb \pr_{12},\Xb \pr_{22})\pr, \ldots, (\Xb \pr_{1n},\Xb \pr_{2n})\pr$ of~$n$ i.i.d.~copies of some~$(d_1+d_2)\!=\!d$-dimensional random vector $(\Xb \pr_1,\Xb \pr_2)\pr$ with Lebesgue-absolutely continuous joint distribution~$\rP\in{\cal P}_d$ and density $f$.  We are interested in the null hypothesis under which $\Xb _1$ and $\Xb _2$, with unspecified marginal distributions~${\rm P}_1\in{\cal P}_{d_1}$ (density~$f_1$) and ${\rm P}_2\in{\cal P}_{d_2}$ (density~$f_2$), respectively, are mutually independent: $f$ then factorizes into $f=f_1f_2$. 

For $k=1,2$ and $i=1,2,\ldots,n$, denote by $R\n_{ki;\pms}$, $\widetilde{R}\n_{ki;\pms}$, and ${\bf S}\n_{ki;\pms}$ the center-outward rank, rescaled rank, and  sign of~$\Xb _{ki}$ computed from~$\Xb _{k1}, \Xb _{k2}, \ldots, \Xb _{kn}$ and the grid $\mathfrak{G}_n$.  Recall that we throughout assume that if~${\bf u}\in \mathfrak{G}_n$ then ${-\bf u}\in\mathfrak{G}_n$.  This implies that $\sum_{i=1}^n{\bf S}\n_{ki;\pms} ={\bf 0}$ and $\sum_{i=1}^n J_k\big(\widetilde{R}\n_{ki;\pms}\big){\bf S}\n_{ki;\pms} = {\bf 0}$ for any {\it score func\-tion}~$J_k: [0, 1) \to \mathbbm{R}$, $k=1,2$.  

Consider now the $d_1\times d_2$ matrices 
\begin{align}
{\tenq{\mathbf W}}\,\!_{\text{\tiny\rm sign}}\n&\coloneqq \frac{1}{n} \sum_{i=1}^n{\bf S}\n_{1i;\pms}{\bf S}\npr_{2i;\pms},\label{tildeW:sign}\\
{\tenq{\mathbf W}}\,\!_{\text{\tiny\rm Spearman}}\n&\coloneqq  \frac{1}{n} \sum_{i=1}^n\widetilde{R}\n_{1i;\pms}\widetilde{R}\n_{2i;\pms}{\bf S}\n_{1i;\pms}{\bf S}\npr_{2i;\pms},\label{tildeW:Spe}\\
{\tenq{\mathbf W}}\,\!_{\text{\tiny\rm Kendall}}\n&\coloneqq  {n \choose 2}^{-1}
\sum_{i<i\pr}
  \text{sign}\Big[
\Big(\widetilde{R}\n_{1i;\pms}{\bf S}\n_{1i;\pms} - \widetilde{R}\n_{1i\pr;\pms}{\bf S}\n_{1i\pr;\pms}\Big)
\Big(\widetilde{R}\n_{2i;\pms}{\bf S}\n_{2i;\pms} - \widetilde{R}\n_{2i\pr;\pms}{\bf S}\n_{2i\pr;\pms}\Big)\pr \Big],\label{tildeW:Ken}
\end{align}
 where  sign$\big[{\bf M}\big]$ stands for the matrix collecting the signs of the entries of a real matrix $\bf M$. More generally, let 
\begin{align}\label{scoreW}
{\tenq{\mathbf W}}\,\!_{J}\n&\coloneqq 
\frac{1}{n} \sum_{i=1}^n
J_1\Big({\widetilde{R}\n_{1i;\pms}}\Big)
J_2\Big({\widetilde{R}\n_{2i;\pms}}\Big)
{\bf S}\n_{1i;\pms}{\bf S}\npr_{2i;\pms}, 
\end{align}
where the {\it score functions} $J_k:[0,1)\to \mathbbm{R}$, $k=1,2$ are continuous and square-integrable, with 
\begin{equation}\label{scorevariance}
0<\lim_{n\to\infty} n^{-1}\sum_{i=1}^{n}J_k^2(r_i) = \int_0^1 J_k^2(u)\rd u \eqqcolon \sigma_{J_k}^2<\infty
\end{equation}
for any $\{r_1,\dots,r_n\}$ such that the uniform distribution over which converges weakly to the uniform distribution over $[0,1]$.
This assumption on score functions will be made throughout the remainder of the paper.
The matrices defined in \eqref{tildeW:sign}--\eqref{scoreW} clearly constitute matrices of cross-covariance measurements based on center-outward ranks and signs (for ${\tenq{\mathbf W}}\,\!_{\text{\tiny\rm sign}}\n$, signs only). 
For~$d_1=1=d_2$, it is easily seen that~${\tenq{\mathbf W}}\,\!_{\text{\tiny\rm sign}}\n$, ${\tenq{\mathbf W}}\,\!_{\text{\tiny\rm Spearman}}\n$, and~${\tenq{\mathbf W}}\,\!_{\text{\tiny\rm Kendall}}\n$, up to scaling constants, reduce to the classical quadrant, Spearman, and Kendall test statistics, while ${\tenq{\mathbf W}}\,\!_{J}\n$ yields a matrix-valued score-based extension of  Spearman's correlation coefficient.

\subsection{Asymptotic representation and asymptotic normality}\label{asreprsec}

Clearly, neither the ranks nor the signs are mutually independent ($n-1$ rank and sign pairs indeed determine the last pair): traditional central-limit theorems thus do not apply.  However, as we show now, each of the rank-based matrices defined in~\eqref{tildeW:sign}--\eqref{scoreW} admits an asymptotic representation in terms of i.i.d.~variables.  More precisely, defi\-ning~${\bf S}_{ki;\pms}$ as~${\bf F}_{k;\pms}({\bf X}_{ki})\big/\big\Vert {\bf F}_{k;\pms}({\bf X}_{ki})\big\Vert$ if ${\bf F}_{k;\pms}({\bf X}_{ki})\neq{\bf 0}$ and ${\bf 0}$ otherwise for $k=1,2$, let 
\begin{align}
{{\mathbf W}}\,\!_{\text{\tiny\rm sign}}\n&\coloneqq \frac{1}{n} \sum_{i=1}^n{\bf S}_{1i;\pms}{\bf S}\pr_{2i;\pms},\label{tildeWas-sign}\\
{{\mathbf W}}\,\!_{\text{\tiny\rm Spearman}}\n&\coloneqq  
\frac{1}{n} \sum_{i=1}^n {\bf F}_{1;\pms}({\bf X}_{1i}) 
 {\bf F}\pr_{2;\pms}({\bf X}_{2i}),\label{tildeWas-Spe}\\
{{\mathbf W}}\,\!_{\text{\tiny\rm Kendall}}\n&\coloneqq  {n \choose 2}^{-1} \sum_{i<i\pr}
  \text{sign}\Big[
\Big({\bf F}_{1;\pms}({\bf X}_{1i}) - {\bf F}_{1;\pms}({\bf X}_{1i\pr}) \Big)
\Big({\bf F}_{2;\pms}({\bf X}_{2i}) - {\bf F}_{2;\pms}({\bf X}_{2i\pr}) \Big)\pr\, \Big],\label{tildeWas-Ken}
\end{align}
 and
\begin{align}\label{scoreWas}
\qquad{{\mathbf W}}\,\!_{J}\n\coloneqq \frac{1}{n} \sum_{i=1}^n
J_1\Big(\big\Vert{\bf F}_{1;\pms}({\bf X}_{1i}) \big\Vert\Big)
J_2\Big(\big\Vert {\bf F}_{2;\pms}({\bf X}_{2i}) \big\Vert\Big) 
{\bf S}_{1i;\pms}{\bf S}\pr_{2i;\pms} .
\end{align}
The following asymptotic representation results then hold under the
null hypothesis of independence (hence, also under contiguous
alternatives).   As usual, let~$\text{\rm vec}({\mathbf A})\coloneqq ({\mathbf a}_1^\prime,\dots,{\mathbf a}_m^\prime)^\prime$ for an $n\times~\!m$ \linebreak matrix ${\mathbf A}=({\mathbf a}_1,\dots,{\mathbf a}_m)$ and   
   write~${\bf X}\n = o_{\text{\rm q.m.}}(c\n)$ when, for some sequence $c\n$  of positive real numbers,  $\vert X\n_j \vert \leq Y\n c\n$ for each component  $X\n_j $ of ${\bf X}\n$, with  $\lim_{n\to\infty} \mathrm{E}[(Y\n)^2]= 0$.

\begin{lemma}\label{lem:hajek} 
Under the null 
hypothesis of independence, as $n$ tends to infinity, it holds that \linebreak $\text{\rm vec}\big({\tenq{\mathbf W}}\,\!_{\text{\tiny\rm sign}}\n \! - {{\mathbf W}}\,\!_{\text{\tiny\rm sign}}\n\big)$, 
$\text{\rm vec}\big({\tenq{\mathbf W}}\,\!_{\text{\tiny\rm Spearman}}\n\! -
{{\mathbf W}}\,\!_{\text{\tiny\rm Spearman}}\n\big)$,
$\text{\rm vec}\big({\tenq{\mathbf W}}\,\!_{\text{\tiny\rm Kendall}}\n\! -
{{\mathbf W}}\,\!_{\text{\tiny\rm Kendall}}\n\big)$, and
$\text{\rm vec}\big({\tenq{\mathbf W}}\,\!_{J}\n - {{\mathbf W}}\,\!_{J}\n\big)$  all\linebreak  are~$\,o_{\text{\rm q.m.}}(n^{-1/2})$. \end{lemma}

  The asymptotic normality of {\rm vec}$({\tenq{\mathbf W}}\,\!_{\text{\tiny\rm sign}}\n)$, {\rm vec}$({\tenq{\mathbf W}}\,\!_{\text{\tiny\rm Spearman}}\n)$, {\rm vec}$({\tenq{\mathbf W}}\,\!_{\text{\tiny\rm Kendall}}\n)$, and {\rm vec}$({\tenq{\mathbf W}}\,\!_{J}\n)$ follows immediately from the asymptotic representation results and the standard central-limit behavior of~{\rm vec}$({{\mathbf W}}\,\!_{\text{\tiny\rm sign}}\n)$, {\rm vec}$({{\mathbf W}}\,\!_{\text{\tiny\rm Spearman}}\n)$, {\rm vec}$({{\mathbf W}}\,\!_{\text{\tiny\rm Kendall}}\n)$, and {\rm vec}$({{\mathbf W}}\,\!_{J}\n)$.

\begin{proposition}\label{prop:asym} 
Under the null hypothesis of independence, as $n$ tends to infinity, it holds that \linebreak $n^{1/2}{\rm vec}({\tenq{\mathbf W}}\,\!_{\text{\tiny\rm sign}}\n),$ 
$n^{1/2}{\rm vec}({\tenq{\mathbf W}}\,\!_{\text{\tiny\rm Spearman}}\n),$ 
$n^{1/2}{\rm vec}({\tenq{\mathbf W}}\,\!_{\text{\tiny\rm Kendall}}\n),$ 
and $n^{1/2}{\rm vec}({\tenq{\mathbf W}}\,\!_{J}\n)$
are asymptotically normal with mean vectors ${\bf 0}_{d_1d_2}$ 
and covariance matrices \vspace{-1mm}
$$\frac{1}{d_1d_2}{\bf I}_{d_1d_2}, \quad 
  \frac{1}{9d_1d_2}{\bf I}_{d_1d_2}, \quad 
  \frac{4}{9}{\bf I}_{d_1d_2}, \quad 
  \text{and} \quad \frac{\sigma^{2}_{J_1}\sigma^{2}_{J_2}}{d_1d_2}{\bf I}_{d_1d_2},
\vspace{-1mm}$$
respectively. 
\end{proposition}

\subsection{Center-outward sign, Spearman, Kendall, and score   tests}\label{testprocsec}

 Associated with   ${\tenq{\mathbf W}}\,\!_{\text{\tiny\rm sign}}\n$, ${\tenq{\mathbf W}}\,\!_{\text{\tiny\rm Spearman}}\n$,  ${\tenq{\mathbf W}}\,\!_{\text{\tiny\rm Kendall}}\n$, and ${\tenq{\mathbf W}}\,\!_{J}\n$ are the quadrant or sign, Spearman, Kendall, and score  test statistics 
\begin{align*}
 {\tenq T}\,\!_{\text{\tiny\rm sign}}\n&\coloneqq  nd_1d_2\big\Vert   {\tenq{\mathbf W}}\,\!_{\text{\tiny\rm sign}}\n\big\Vert^2_{\mathrm F},\quad &
  {\tenq T}\,\!_{\text{\tiny\rm Spearman}}\n&\coloneqq  9nd_1d_2\big\Vert   {\tenq{\mathbf W}}\,\!_{\text{\tiny\rm Spearman}}\n\big\Vert^2_{\mathrm F}, \quad \\
 {\tenq T}\,\!_{\text{\tiny\rm Kendall}}\n&\coloneqq  \frac{9n}{4}\big\Vert
   {\tenq{\mathbf W}}\,\!_{\text{\tiny\rm Kendall}}\n\big\Vert^2_{\mathrm
   F},  \qquad\qquad \text{and} &\hspace{-5mm}
  {\tenq T}\,\!_{J}\n&\coloneqq 
\frac{nd_1d_2}{\sigma^{2}_{J_1}\sigma^{2}_{J_2}}
\big\Vert   {\tenq{\mathbf W}}\,\!_J\n\big\Vert^2_{\mathrm F},   
\end{align*}
respectively, where $\Vert {\bf M}\Vert _{\mathrm F}$ denotes the Frobenius
norm of a matrix $\bf M$, and $\sigma^{2}_{J_1}$, $\sigma^{2}_{J_2}$ are defined as in~\eqref{scorevariance}. 
{
In particular, considering the 
van der Waerden\vspace{-1mm} score functions~$J^{\text{\tiny{\rm vdW}}}_k(u)\!\coloneqq \!\big(F^{-1}_{\chi^2_{d_k}}\!(u)\big)^{\! 1/2}\!$, where~$F_{\chi^2_d}$ stands for the $\chi^2_d$ distribution function, the van der Waerden test statistic is~${\tenq T}\,\!_{\text{\rm{\tiny vdW}}}\n=\,\,{\tenq T}\,\!_{J^{\text{\tiny{\rm vdW}}}}\n =n\big\Vert   {\tenq{\mathbf W}}\,\!_{J^{\text{\tiny{\rm vdW}}}}\n\big\Vert^2_{\mathrm F}$.
\vspace{1mm}}

In view of the asymptotic normality results in Proposition~\ref{prop:asym}, the tests 
 ${\psi}\,\!_{\text{\tiny\rm sign}}\n$, ${\psi}\,\!_{\text{\tiny\rm Spearman}}\n$, ${\psi}\,\!_{\text{\tiny\rm Kendall}}\n$, {${\psi}\,\!_{\text{\tiny\rm vdW}}\n$}, and~${\psi}\,\!_{J}\n$ rejecting the null hypothesis of independence whenever ${\tenq T}\,\!_{\text{\tiny\rm sign}}\n$, ${\tenq T}\,\!_{\text{\tiny\rm Spearman}}\n$, ${\tenq T}\,\!_{\text{\tiny\rm Kendall}}\n$, {${\tenq T}\,\!_{\text{\tiny\rm vdW}}\n$} or ${\tenq T}\,\!_{J}\n$, respectively, exceed the~$(1-~\!\alpha)$-quan\-tile~$\chi^2_{d_1d_2;1-\alpha}$ of a chi-square distribution with~$d_1d_2$ degrees of freedom have asymptotic\footnote{In view of the  distribution-freeness of their test statistics, the convergence to $\alpha$ of the actual size of these tests is uniform over the possible distributions~${\rm
P}_1\in{\cal P}_{d_1}$  and ${\rm
P}_2\in{\cal P}_{d_2}$. } level $\alpha$. These tests are, however, strictly distribution-free, and exact critical values can be computed or simulated as well. The tests based on~${\tenq T}\,\!_{\text{\tiny\rm sign}}\n$, ${\tenq T}\,\!_{\text{\tiny\rm Spearman}}\n$, ${\tenq T}\,\!_{\text{\tiny\rm Kendall}}\n$, and {${\tenq T}\,\!_{\text{\tiny\rm vdW}}\n$} are multivariate extensions of the traditional quadrant, Spearman, Kendall and {van der Waerden} tests, respectively,
 to which they reduce for~$d_1~\!=~\!1=~\!d_2$.

\section{Local asymptotic powers}\label{sec:4}

While there is only one way for two random vectors ${\bf X}_1$
and~${\bf X}_2$ to be independent, their mutual dependence can take
many forms. The classical benchmark, in testing for bivariate
independence, is a ``local'' form of an independent component analysis model that goes back to \citet{MR79384} and \citet{MR0229351}. Multivariate elliptical extensions of such alternatives have been considered also by \citet{MR1467849}, \citet{MR1965367}, \citet{MR2462206}, and  \citet{MR3544291}.

\subsection{Generalized Konijn families}\label{Konijnsec}

Let ${\bf X} ^*=({\bf X} ^{*\prime}_1, {\bf X} ^{*\prime}_2)\pr$,
where ${\bf X} ^*_1$ and ${\bf X} ^*_2$  are mutually independent
random vectors, with absolutely continuous distributions $\rP_{1}$
over $\mathbbm{R}^{d_1}$ and~$\rP_{2}$ over~$\mathbbm{R}^{d_2}$ and
densities $f_1$ and $f_2$, respectively.  Then~${\bf X} ^*$ has density  $f=f_1f_2$ over $\mathbbm{R}^{d}$. Consider 
\begin{align}\label{Konijn}
{\bf X}=\left(\!
\begin{array}{c}
{\bf X}_1\\  {\bf X}_2
\end{array}
\!\right)
\coloneqq {\bf M}_\delta\left(\!
\begin{array}{c}
{\bf X} ^*_1\\  {\bf X} ^*_2
\end{array}
\!\right) 
\coloneqq 
\left(\!
\begin{array}{cc}
(1-{\delta}
){\bf I}_{d_1}&{\delta}{\bf M_1}\\
{\delta}{\bf M}_2&(1-{\delta}){\bf I}_{d_2}
\end{array}
\!\right)
\!\left(\!
\begin{array}{c}
{\bf X} ^*_1\\  {\bf X} ^*_2
\end{array}
\!\right)
\end{align}
where ${\delta}\!\in\!\mathbbm{R}$ and ${\bf M}_1\!\in\!\mathbbm{R}^{d_1\times d_2}$, ${\bf M}_2\!\in\!\mathbbm{R}^{d_2\times d_1}$ are nonzero matrices.  For given $\rP_{1}$, $\rP_{2}$, ${\bf M}_1$, ${\bf M}_2$, and~$\delta$, the distribution~$ {\rm P}^{\bf X}_{ \delta}({\rm P}_1, {\rm P}_2; {\bf M}_1, {\bf M}_2)$ of~${\bf X}$ belongs to the one-parameter family 
\[
{\cal P}^{\bf X}_{{\rm P}_1, {\rm P}_2; {\bf M}_1, {\bf M}_2}\coloneqq \left\{\rP^{\bf X}_{\delta} ({\rm P}_1, {\rm P}_2; {\bf M}_1, {\bf M}_2)\vert\, {\delta}\in\mathbbm{R}\right\};
\]
call it a {\it generalized Konijn family}. Independence between ${\bf X}_1$ and ${\bf X}_2$, in such families, holds iff $\delta = 0$. 
  
Sequences of the form~${\rm P}^{\bf X}_{{n^{-1/2}\tau}} ({\rm P}_1, {\rm P}_2; {\bf M}_1, {\bf M}_2)\vspace{1mm}$ with $\tau\neq 0$, as we shall see, define contiguous alternatives to
the null hypothesis of independence in a sample of size~$n$.  More precisely, consider the following regularity assumptions on ${\rm P}_1$ and ${\rm P}_2$ (on their densities $f_1$ and $f_2$). 
\begin{assumption}\label{asp:K}
(K1) The densities $f_1$ and $f_2$ are such that
\begin{align*}
\int_{\mathbbm{R}^{d_k}}\!{\bf x}f_k({\bf x})\rd {\bf x}={\bf 0}\quad\text{and}\quad
0<\int_{\mathbbm{R}^{d_k}}\!{\bf x}{\bf x}\pr f_k({\bf x})\rd {\bf x}\eqqcolon {\boldsymbol{\Sigma}}_{k}<\infty,\quad k=1,2,
\end{align*}
where $0< {\boldsymbol{\Sigma}}_{k}<\infty,$ means the covariance matrix ${\boldsymbol\Sigma}_k$ is finite and positive definite.
\begin{compactenum}
\item[(K2)] The functions ${\bf x}_k\mapsto (f_k({\bf x}_k))^{1/2}$, $k=1,2$
admit quadratic mean partial derivatives\footnote{Existence of quadratic mean partial derivatives  is equivalent to quadratic mean differentiability; this was shown in \citet{MR307329} and independently rediscovered by  \citet[Lemma~2.1]{MR1364260}.}\vspace{-1mm} 
$$D_{\ell}[(f_k)^{1/2}], \quad \ell =1,\ldots,d_k, \ k=1,2.\vspace{-1mm}$$
\item[(K3)] The score ${\boldsymbol\varphi}\mkern-2mu\coloneqq\mkern-2mu  \left({\boldsymbol\varphi}_1\pr ,{\boldsymbol\varphi}_2\pr\right)\pr\mkern-2mu\coloneqq\mkern-2mu \left(\varphi_{1;1},\ldots,\varphi_{1;d_1},\varphi_{2;1},\ldots,\varphi_{2;d_2}\right)\pr$ with $\varphi_{k;\ell}\mkern-2mu\coloneqq\mkern-2mu -2D_\ell[(f_k)^{1/2}]/(f_k)^{1/2}$ is such that $0<~\int_{{\mathbbm R}^{d_k} }\big( \varphi_{k;\ell}({\bf x})\big)^{2}f_k({\bf x})\rd {\bf x}<\infty$, $k=1,2$, $ \ell =1,\ldots,d_k$,\footnote{Integration by parts yields $\int_{{\mathbbm R}^{d_k} }             {\boldsymbol\varphi}_k({\bf x})    f_k({\bf x})\rd {\bf x} ={\bf 0}$,
$\int_{{\mathbbm R}^{d_k} } {\bf x}\pr {\boldsymbol\varphi}_k({\bf x})    f_k({\bf x})\rd {\bf x} =d_k$, and 
$\int_{{\mathbbm R}^{d_k} } {\bf x}    {\boldsymbol\varphi}_k({\bf x})\pr f_k({\bf x})\rd {\bf x} =~{\bf I}_{d_k}$, \linebreak 
for~$k=1,2$; see also \citet[page~555]{MR1364260}.} and
\begin{align*}
\mathcal{J}_k
\coloneqq {\text{\rm Var}}\left({\bf X} ^{*\prime}_k {\boldsymbol\varphi}_k({\bf X} ^*_k)\right)=
  \int_{{\mathbbm R}^{d_k} } \left({\bf x}\pr {\boldsymbol\varphi}_k({\bf x}) -d_k\right)^2
  f_k({\bf x})\rd {\bf x},\quad k=1,2\vspace{-2mm}
\end{align*}
is finite and strictly positive.\vspace{-2mm}
\end{compactenum}
\end{assumption}
\noindent It should be stressed that these regularity assumptions are not to be imposed on the observations in order for our rank-based tests to be valid; they only are required for deriving local power and ARE results. \smallskip

The following local asymptotic normality (LAN) property then holds in the vicinity of ${\delta}=0$. \pagebreak

\begin{lemma}\label{lem:lan}  Denote by ${\rm P}\n_{\delta}$\!, $n\in\mathbbm{N}$,  the distribution of the triangular array ${\bf X}\n\coloneqq ({\bf X}_1\n,\ldots,{\bf X}_n\n)$  of~$n$ independent copies of\, ${\bf X}=({\bf X}^\prime_1, {\bf X}_2^\prime)^\prime\sim\rP^{\bf X}_{\delta} ({\rm P}_1, {\rm P}_2; {\bf M}_1, {\bf M}_2)\in {\cal P}^{\bf X}_{{\rm P}_1, {\rm P}_2; {\bf M}_1, {\bf M}_2}$, where ${\rm P}_1$ and ${\rm P}_2$ \linebreak satisfy Assumption~4.1. Then, the fa\-mily~$\big\{{\rm P}\n_{\delta}\vert\, \delta\in{\mathbb R}\big\}$ is   Locally Asymptotically Normal (LAN) at~${\delta}=~\!0$, with root-$n$ contiguity rate, central sequence\vspace{-1mm}
\begin{align*}
\Delta \n({\bf X}\n)
\coloneqq n^{-1/2}
\sum_{i=1}^n
\left[
{\bf X}^{(n)\prime}_{1i}
{\bf M}_2\pr{\boldsymbol\varphi}_2({\bf X}\n_{2i})
\right. &+ 
{\bf X}^{(n)\prime}_{2i}{\bf M}_1\pr{\boldsymbol\varphi}_1({\bf X}\n_{1i}) 
\\
&\left. -\Big({\bf X}^{(n)\prime}_{1i}{\boldsymbol\varphi}_1({\bf X}\n_{1i}) -d_1
\Big)
-\Big({\bf X}^{(n)\prime}_{2i}{\boldsymbol\varphi}_2({\bf X}_{2i}\n) -d_2
\Big)
\right],
\end{align*}
and Fisher information 
\begin{equation*}
\gamma ^2\coloneqq {\cal J}_1 + {\cal J}_2 
 + \text{\rm vec} \pr\! \left({\boldsymbol{\Sigma}}_1\right) \text{\rm vec}\! \left({\bf M}_2\pr{\bm{\mathcal I}}_2{\bf M}_2\right) 
 + \text{\rm vec} \pr\! \left({\boldsymbol{\Sigma}}_2\right) \text{\rm vec}\! \left({\bf M}_1\pr{\bm{\mathcal I}}_1{\bf M}_1\right) 
 + \text{\rm tr}({\bf M}_1{\bf M}_2) 
 + \text{\rm tr}({\bf M}_2{\bf M}_1)
\end{equation*}
where 
${\bm{\mathcal I}}_k\coloneqq \int_{{\mathbbm R}^{d_k} }{\boldsymbol\varphi}_k({\bf x}){\boldsymbol\varphi}_k\pr({\bf x}) f_k({\bf x})\rd {\bf x}$, $k=1,2$. 
In other words, under ${\rm P}^{(n)}_0$, as $n\to\infty$, 
\begin{equation}
\Lambda\n({\bf X}\n)\coloneqq \log\frac{\rd\rP^{(n)}_{n^{-1/2}\tau}}{\rd\rP^{(n)}_0}({\bf X}\n) 
=\tau\Delta \n({\bf X}\n) -\frac{1}{2}\tau^2\gamma ^2 + o_{\rP}(1)
\label{LANKon}
\end{equation}
and $\Delta \n({\bf X}\n)$ is asymptotically normal, with mean zero and variance $\gamma ^2$.
 \end{lemma}

If  ${\rm P}_1$  and ${\rm P}_2$ are elliptical with mean $\boldsymbol 0$, scatter matrices   ${\boldsymbol \Sigma}_k$, and     radial densities $\phi_k$,  $k=~\!1,2$, viz.
$$f_k({\bf x}_k)\propto ({\rm det}({\boldsymbol{\Sigma}}_k))^{-1/2} 
\phi_k\Big(\sqrt{{\bf x}_k^{\prime}{\boldsymbol{\Sigma}}_k^{-1}{\bf x}_k}\,\Big),\quad k=1,2,$$
write~${\rm P}_\delta^{\,\text{\rm{\tiny ell}}}(\phi_1,\phi_2, {\boldsymbol \Sigma}_1, {\boldsymbol \Sigma}_2; {\bf M}_1, {\bf M}_2)$ and ${\cal P}^{\,\text{\rm{\tiny ell}}}_{\phi_1,\phi_2, {\boldsymbol \Sigma}_1, {\boldsymbol \Sigma}_2; {\bf M}_1, {\bf M}_2}\vspace{1mm}$, respectively,  instead of~${\rm P}_\delta^{\bf X}({\rm P}_1,{\rm P}_2; {\bf M}_1, {\bf M}_2)$ and~${\cal P}^{\bf X}_{{\rm P}_1,{\rm P}_2; {\bf M}_1, {\bf M}_2}\vspace{1mm}$.  The Konijn alternatives and Konijn families considered in  \citet{MR1467849}, \citet{MR1965367}, and \citet{MR2462206} are particular cases, of the form ${\rm P}_\delta^{\,\text{\rm{\tiny ell}}}(\phi_1,\phi_2, {\bf I}_{d_1}, {\bf I}_{d_2}; {\bf M}, {\bf M}\pr)$ and~${\cal P}^{\text{\rm{\tiny ell}}}_{\phi_1,\phi_2, {\bf I}_{d_1}, {\bf I}_{d_2}; {\bf M}, {\bf M}\pr}$: call them {\it elliptical} Konijn alternatives and {\it elliptical}  Konijn families, respectively.

\subsection{Limiting distributions and Pitman efficiencies}

For the univariate two-sample location problem, \citet{MR100322} and \citet{MR79383} established their celebrated results on the asymptotic relative efficiency (ARE) of traditional normal score (van der Waerden) and Wilcoxon rank tests with respect to Student's Gaussian procedure. These results were extended by \cite{MR1312324} and \citet{MR1833863} to ARMA time-series models, by \citet{MR1926170} to Mahalanobis ranks-and-signs-based location tests under multivariate elliptical distributions, by \citet{MR1963662} to VARMA time-series models with elliptical innovations, by  \citet{MR2160622} for the shape parameter of elliptical distributions. 

Chernoff--Savage and Hodges--Lehmann bounds for the AREs (with respect to Hotelling) of measure-transportation-based center-outward rank and sign tests were first obtained by  \citet{deb2021efficiency} in the context of two-sample location models and for the subclasses of elliptical and independent component distributions.  In this section, we similarly establish Chernoff--Savage and Hodges--Lehmann results for the AREs, relative to Wilks' test, of our center-outward van der Waerden and Spearman tests, respectively. The two-sample location model here is replaced with an elliptical Konijn model satisfying  (in order for ARE values to make sense) Assumptions~4.1.  

To this end, we first derive the limiting distributions of~${\tenq T}\,\!_{J}\n$ (of which~${\tenq T}\,\!_{\text{\rm{\tiny sign}}}\n$,  ${\tenq T}\,\!_{\text{\rm{\tiny Spearman}}}\n$, and ${\tenq T}\,\!_{\text{\rm{\tiny vdW}}}\n$ are particular cases) and ${\tenq T}\,\!_{\text{\tiny\rm Kendall}}\n$ under local sequences of {(not necessarily elliptical)} Konijn alternatives of the form~${\rm P}^{\bf X}_{n^{-1/2}\tau}({\rm P}_1, {\rm P}_2; {\bf M}_1, {\bf M}_2)$. 
\begin{theorem}\label{thm:J}
Let $\rP_{1}$ and $\rP_{2}$ satisfy Assumption~\ref{asp:K}.  If the observations are $n$ independent copies of~${\bf X}\sim {\rm P}^{\bf X}_{n^{-1/2}\tau}({\rm P}_1, {\rm P}_2; {\bf M}_1, {\bf M}_2)$, 
\begin{compactenum}
\item[(i)] the limiting distribution as $n\to\infty$ of the test statistic ${\tenq T}\,\!_{J}\n$ is noncentral chi-square with $d_1d_2$ degrees of freedom and noncentrality parameter 
\begin{equation}\label{still}
C_{J}(\tau)\Big\Vert {\rm E}_{H_0}\Big[{\bf J}_1({\bf F}_{1;\pms} ({\bf X}_1))
{\bf R} {\bf J}_2({\bf F}_{2;\pms} ({\bf X}_2))\pr\Big]\Big\Vert^2_{\mathrm F},
\end{equation}
where $C_{J}(\tau)\coloneqq \frac{\tau^2 d_1d_2}{\sigma^{2}_{J_1}\sigma^{2}_{J_2}}$,  ${\bf R}\coloneqq {\bf X}_{1}\pr{\bf M}_2\pr{\boldsymbol\varphi}_2({\bf X}_{2}) + 
{\bf X}_{2}\pr{\bf M}_1\pr{\boldsymbol\varphi}_1({\bf X}_{1})$, 
 $\displaystyle{{\bf J}_k({\bf u})\coloneqq  J_k(\Vert{\bf u}\Vert)\frac{{\bf u}}{\Vert{\bf u}\Vert}{\bf 1}_{[\Vert{\bf u}\Vert\neq 0]}}$ for ${\bf u}\in{\mathbb{S}_d}$,
 ${\rm E}_{H_0}$ stands for expectations under the null ($\delta = 0$),
 and recall that $\Vert{\bf M}\Vert _{\rm F}$ denotes the Frobenius norm~of~$\bf M$;
\item[(ii)] the limiting distribution of the Kendall test statistic ${\tenq T}\,\!_{\text{\tiny\rm Kendall}}\n$ is noncentral chi-square with $d_1d_2$ degrees of freedom and noncentrality parameter
\begin{equation*}
{9\tau^2}\Big\Vert   {\rm E}_{H_0}\Big[{\bf F}^{\scriptscriptstyle\square}_{1;\pms} ({\bf X}_1)
{\bf R} {\bf F}^{\scriptscriptstyle\square}_{2;\pms} ({\bf X}_2)\pr\Big]\Big\Vert^2_{\mathrm F},
\end{equation*}
where, denoting by $F_{k\ell}$ the marginal cumulative distribution function of~$\big({\bf F}_{k;\pms}({\bf X}_{k}) \big)_{\ell}$ and
$$\big({\bf F}^{\scriptscriptstyle\square}_{k;\pms} ({\bf X}_k)\big)_{\ell}\coloneqq 2F_{k\ell}\big(\big({\bf F}_{k;\pms}({\bf X}_{k})\big)_{\ell}\big)-1,\qquad k=1,2,\ \ \ell=1,\ldots,d_k.$$  
\end{compactenum}
\end{theorem}
\noindent The results for ${\tenq T}\,\!_{\text{\rm{\tiny sign}}}\n$, ${\tenq T}\,\!_{\text{\rm{\tiny Spearman}}}\n$, and ${\tenq T}\,\!_{\text{\rm{\tiny vdW}}}\n$ follow as particular cases of (i) with the constant  $C_J(\tau)$ taking values $\tau^2 d_1d_2$, $9\tau^2 d_1d_2$, and $\tau^2$, respectively.\smallskip

Turning to Wilks' test, the test statistic takes the form 
\[
T\n_{\text{\tiny{\rm Wilks}}}\coloneqq n\log V\n\quad\text{ with }\quad 
V\n\coloneqq \frac{{\rm det}({\bf S}_1\n){\rm det}({\bf S}_2\n)}{{\rm det}({\bf S}\n)}, 
\]
where ${\rm det}({\bf M})$ denotes the determinant of a square matrix $\bf M$, ${\bf S}_k\n$ is the $d_k\times d_k$ sample covariance matrix  of~$\Xb_{k1}, \ldots, \Xb_{kn}$, $k=1,2$, and ${\bf S}\n$  the sample covariance matrix of $(\Xb \pr_{11},\Xb \pr_{21})\pr, \ldots, (\Xb \pr_{1n},\Xb \pr_{2n})\pr$, which is $(d_1+d_2)\times (d_1+d_2)$. Under the null hypothesis of independence, $T\n_{\text{\tiny{\rm Wilks}}}$ is {asymptotically} chi-square with~$d_1d_2$ degrees of freedom as soon as ${\bf X}_1$ and~${\bf X}_2$ have finite variances  {(\citealp[pages~187--188]{MR2462206}; \citealp[Section~5.1]{MR3544291})}. Wilks' test $\psi\n_{\text{\tiny{\rm Wilks}}}$, thus,  rejects the null hypothesis at asymptotic level $\alpha$ whenever~$T\n_{\text{\tiny{\rm Wilks}}}$ exceeds the corresponding  chi-square quantile of order $(1-\alpha)$. Under alternatives of the form~${\bf X}\sim {\rm P}_{n^{-1/2}\tau}^{\,\text{\rm{\tiny ell}}}(\phi_1,\phi_2,{\boldsymbol\Sigma}_1, {\boldsymbol\Sigma}_2; {\bf M}_1, {\bf M}_2)\vspace{1mm}$, with $f_1$ and $f_2$ satisfying Assumption~\ref{asp:K}, the limiting distribution of Wilks' statistic is noncentral chi-square, still with~$d_1d_2$ degrees of freedom, with noncentrality parameter
\begin{equation}\label{still2}
\tau^2 \Big\Vert 
{\boldsymbol{\Sigma}}_{1}^{ 1/2}{\bf M}_2\pr{\boldsymbol{\Sigma}}_2^{-1/2} + 
{\boldsymbol{\Sigma}}_{1}^{-1/2}{\bf M}_1   {\boldsymbol{\Sigma}}_2^{ 1/2} \Big\Vert^2_{\mathrm F}
\end{equation}
{reducing, under the elliptical Konijn alternative $\vspace{1mm}{\rm P}_{n^{-1/2}\tau}^{\,\text{\rm{\tiny ell}}}(\phi_1,\phi_2, {\bf I}_{d_1}, {\bf I}_{d_2}; {\bf M}, {\bf M}\pr)$, to $4\tau^2 \Vert 
{\bf M}\Vert^2_{\mathrm F}$} (see, e.g., page 919 of \citet{MR2201019}). 
We are now ready to derive asymptotic relative efficiency results for our center-outward rank tests relative to Wilks' pseudo-Gaussian one. 

\begin{proposition}\label{prop:Pitman}
Let $\rP_{1}$ and $\rP_{2}$ satisfying Assumption~\ref{asp:K} be elliptically symmetric with radial densities $\phi_k$ and covariance matrices~${\boldsymbol\Sigma}_k$, $k=1,2$.  
 Then, 
 the Pitman asymptotic relative efficiency (ARE) in~${\cal P}^{\,\text{\rm{\tiny ell}}}_{\phi_1,\phi_2, {\boldsymbol \Sigma}_1, {\boldsymbol \Sigma}_2; {\bf M}_1, {\bf M}_2}\vspace{1mm}$ 
 of the center-outward test $\psi\n_{J}$ based on  the score functions~$J_k$, $k=1,2$
 relative to Wilks' test $\psi\n_{\text{\tiny{\rm Wilks}}}$ is
\begin{align}\label{AREell}
{\rm ARE}({\psi}\,\!_{J}\n, \psi\n_{\text{\tiny{\rm Wilks}}})
=\frac{\Big\Vert 
D_1C_2{\boldsymbol{\Sigma}}_{1}^{ 1/2}{\bf M}_2\pr{\boldsymbol{\Sigma}}_2^{-1/2} + 
D_2C_1{\boldsymbol{\Sigma}}_{1}^{-1/2}{\bf M}_1   {\boldsymbol{\Sigma}}_2^{ 1/2} \Big\Vert^2_{\mathrm F}}{d_1d_2\sigma^{2}_{J_1}\sigma^{2}_{J_2}\Big\Vert 
{\boldsymbol{\Sigma}}_{1}^{ 1/2}{\bf M}_2\pr{\boldsymbol{\Sigma}}_2^{-1/2} + 
{\boldsymbol{\Sigma}}_{1}^{-1/2}{\bf M}_1   {\boldsymbol{\Sigma}}_2^{ 1/2} \Big\Vert^2_{\mathrm F}},
\end{align}
where, for $k=1,2$,
\[
C_k= C_k(J_k,\phi_k)\coloneqq {\rm E}[J_k(U)\rho_k(\widetilde F_k^{-1}(U))]\quad\text{and}\quad
D_k=D_k(J_k,\phi_k)\coloneqq {\rm E}[J_k(U) \widetilde F_k^{-1}(U))]
\]
\noindent with $\rho_k\coloneqq \frac{-\,\phi_k^{\prime}}{\phi_k}$,   
$\widetilde F_k$ the cumulative distribution function of $\Vert  {\boldsymbol{\Sigma}}_{k}^{-1/2}{\bf X}_k\Vert$, 
and~$U$ a random variable uniformly distributed over $(0,1)$. 
  
If, moreover, ${\boldsymbol{\Sigma}}_{1} {\bf M}_2\pr
={\bf M}_1 {\boldsymbol{\Sigma}}_{2}$---hence, in particular, 
under the family ${\cal P}^{\,\text{\rm{\tiny ell}}}_{\phi_1,\phi_2, {\bf I}_{d_1}, {\bf I}_{d_2}; {\bf M}, {\bf M}\pr}\vspace{1mm}$, we have
\begin{compactenum}
   \item[(i)]${\rm ARE}({\psi}\,\!_{{\text{\tiny{\rm vdW}}}}\n, \psi_{\text{\tiny{\rm Wilks}}}\n)\ge 1,$ 
where ${\psi}\,\!_{{\text{\tiny{\rm vdW}}}}\n= \psi_{J^{\text{\tiny{\rm vdW}}}}\n$ with van der Waerden score functions $J^{\text{\tiny{\rm vdW}}}_k,~k=1,2$
and equality under Gaussian ${\rm P}_1$ and ${\rm P}_2$ only;
   \item[(ii)]
$
 {\rm ARE}({\psi}\,\!_{{\text{\tiny{\rm Spearman}}}}\n, \psi\n_{\text{\tiny{\rm Wilks}}})
 \ge \Omega(d_1,d_2)
 \ge {9}/{16},
$
where ${\psi}\,\!_{{\text{\tiny{\rm Spearman}}}}\n= \psi\n_{J^{\text{\tiny{\rm W}}}}$ with the
Wilcoxon score functions $J^{\text{\tiny{\rm W}}}_k(u)\coloneqq u$,\vspace{-2mm} 
\[\Omega(d_1,d_2)\coloneqq  {9(2c_{d_1}^2+d_1-1)^2(2c_{d_2}^2+d_2-1)^2}/{1024 d_1d_2c_{d_1}^2c_{d_2}^2}, \vspace{-2mm}\]
\[c_d\coloneqq \inf\big\{ x>0 \ \big\vert\  \big(\sqrt{x} B_{\sqrt{2d-1} / 2}(x)\big)\pr = 0\big\},\quad\text{where}\quad
B_{a}(x)\coloneqq \sum_{m=0}^{\infty}{\frac {(-1)^{m}}{m!\Gamma (m+a+1)}}{\left({\frac {x}{2}}\right)}^{2m+a}. 
\]
\end{compactenum}
 \end{proposition}

Note that the ARE values in~\eqref{AREell} depend, via ${\boldsymbol{\Sigma}}_{1}$ and ${\boldsymbol{\Sigma}}_{2}$, on the underlying covariance structure   of~${\rm P}_0^{\,\text{\rm{\tiny ell}}}(\phi_1,\phi_2, {\boldsymbol \Sigma}_1, {\boldsymbol \Sigma}_2; {\bf M}_1, {\bf M}_2)$; 
 a similar fact was already pointed out by~\citet{MR2691505}.  
Most authors  (e.g. \citet{MR2691505}, \citet{MR1467849},
\citet{MR1965367,MR2088309}, \citet{MR2201019}, \citet{MR2462206} and \citet{deb2021efficiency}), { therefore,   restrict themselves  to} the particular  case of  elliptical Konijn fami\-lies~${\cal P}^{\,\text{\rm{\tiny ell}}}_{\phi_1,\phi_2, {\bf I}_{d_1}, {\bf I}_{d_2}; {\bf M}, {\bf M}\pr}\vspace{0mm}$.  
In \eqref{AREell}, we more generally 
 provide explicit results for the broader class of ${\cal P}^{\,\text{\rm{\tiny ell}}}_{\phi_1,\phi_2, {\boldsymbol\Sigma}_1, {\boldsymbol\Sigma}_2; {\bf M}_1, {\bf M}_2}\vspace{0mm}$ families with arbitrary covariances~${\boldsymbol{\Sigma}}_{1}$, ${\boldsymbol{\Sigma}}_{2}$ and arbitrary matrices ${\bf M}_1$ and~${\bf M}_2$. 
 
Claim (i) entails  the Pitman non-admissibility,  {within the class ${\displaystyle \bigcup}_{\phi_1, \phi_2, {\bf I}_{d_1}, {\bf I}_{d_2}, {\bf M}}
{\cal P}^{\text{\rm{\tiny ell}}}_{\phi_1, \phi_2, {\bf I}_{d_1}, {\bf I}_{d_2}; {\bf M}, {\bf M}\pr}$  
of all elliptical Konijn families satisfying Assumption~4.1,    of Wilks' test, 
which is uniformly dominated by our center-outward   van der Waerden test.  
This is comparable with Theorem~4.1 in \citet{deb2021efficiency}, which deals with two-sample location tests. 
The proof relies on Proposition~1 in \citet{MR2462206} (see also~\citet{MR2160622}).

Claim (ii) is a multivariate extension of \citet{MR79383}'s univariate result. 
The infimum  9/16$\,=\,$0.5625 of~$\Omega(d_1,d_2)$ is achieved for~$d_1,d_2\to\infty$.  Table~{S.2} in the supplement gives numerical values of~$\Omega(d_1,d_2)$ 
for $d_1,d_2\le 10$.

\section{Conclusion}\label{sec:6}

Optimal transport provides a new approach to rank-based statistical inference in dimension $d\geq 2$.  The new multivariate ranks retain many of the favorable properties one is used to from the classical univariate ranks.  Here, we demonstrate how the new multivariate ranks can be used for a definition of multivariate versions of popular rank-based correlation coefficients such as quadrant (sign) correlation, Kendall's tau, or Spearman's rho.  We show how the new multivariate rank correlations yield fully distribution-free, yet powerful and computationally efficient tests of vector independence for which we provide explicit local asymptotic powers.  We also show that the van der Waerden (normal score) version of our tests enjoy, relative to Wilks' classical pseudo-Gaussian procedure and, under an elliptical generalization of the traditional Konijn alternatives, the Chernoff--Savage property of their classical bivariate counterparts---which makes Wilks' test non-admissible in the context.  These results are not specific to our choice of the uniform over the unit ball as the reference measure in our optimal transport approach; transports to the Lebesgue uniform over the unit cube (as in \citet{MR3611491})  or ``direct'' transports to, e.g., the spherical Gaussian (as in \citet{deb2021efficiency}) would yield the same local asymptotic powers and asymptotic relative efficiencies.

\section{Proofs}\label{sec:7}

\noindent {\bf Proof of Lemma~\ref{lem:hajek}.} We only need to prove the result for {\rm vec}$({\tenq{\mathbf W}}\,\!_{J}\n)$ (since~{\rm vec}$({\tenq{\mathbf W}}\,\!_{\text{\tiny\rm sign}}\n)$, {\rm vec}$({\tenq{\mathbf W}}\,\!_{\text{\tiny\rm Spearman}}\n)$, and~{\rm vec}$({\tenq{\mathbf W}}\,\!_{\text{\tiny\rm vdW}}\n)$ are particular cases) and~{\rm vec}$({\tenq{\mathbf W}}\,\!_{\text{\tiny\rm Kendall}}\n)$. \smallskip 

\noindent  (a) Starting with  {\rm vec}${\tenq{\mathbf W}}\,\!_{J}\n$,
 letting 
\begin{equation*}
{\bf Y}\n_{ki}\coloneqq J_k\Big(\widetilde{R}\n_{ki;\pms}\Big){\bf S}\n_{ki;\pms} 
\quad \text{and} \quad 
{\bf Y}_{ki}\coloneqq J_k\Big(\big\Vert{\bf F}_{k;\pms}({\bf X}_{ki}) \big\Vert\Big){\bf S}_{ki;\pms}, \quad k=1,2,~i=1,\dots,n,
\end{equation*}
rewrite $n^{1/2}\big( {\tenq{\mathbf W}}\,\!_{J}\n
-{{\mathbf W}}\,\!_{J}\n \big)_{\ell_1 \ell_2}$, $\ell_1=1,\dots,d_1$, $\ell_2=1,\dots,d_2$ as
\begin{align*}
&n^{-1/2}\sum_{i=1}^n\Big(
 \big({\bf Y}\n_{1i}\big)_{\ell_1}\big({\bf Y}\n_{2i}\big)_{\ell_2}
-\big({\bf Y}  _{1i}\big)_{\ell_1}\big({\bf Y}  _{2i}\big)_{\ell_2}\Big).
\end{align*}
It suffices to show that
\begin{align}
{\rm E}\Big[\Big\{n^{-1/2}\sum_{i=1}^n\big({\bf Y}_{1i}\big)_{\ell_1}\big({\bf Y}_{2i}\big)_{\ell_2}\Big\}^2\Big],
\qquad 
{\rm E}\Big[\Big\{n^{-1/2}\sum_{i=1}^n\big({\bf Y}\n_{1i}\big)_{\ell_1}\big({\bf Y}\n_{2i}\big)_{\ell_2}\Big\}^2\Big],\label{two}
\end{align}
\text{and}
\begin{align}
{\rm E}\Big[\Big\{n^{-1/2}\sum_{i=1}^n\big({\bf Y}\n_{1i}\big)_{\ell_1}\big({\bf Y}\n_{2i}\big)_{\ell_2}\Big\}\Big\{
n^{-1/2}\sum_{i=1}^n\big({\bf Y}_{1i}\big)_{\ell_1}\big({\bf Y}_{2i}\big)_{\ell_2}\Big\}\Big]\label{three}
\end{align}
tend to the same limit as $n$ tends to infinity. 

First, consider the left-hand quantity in~\eqref{two}. Due to the independence between $\big({\bf Y}_{1i}\big)_{\ell_1}$ and $\big({\bf Y}_{2i}\big)_{\ell_2}$ and the independence between $\big({\bf Y}_{1i}\big)_{\ell_1}$ and $\big({\bf Y}_{2i\pr}\big)_{\ell_2}$, $i\neq i\pr$, $\ell_1=1,\dots,d_1$, $\ell_2=1,\dots,d_2$, we have 
\begin{align*}
{\rm E}\Big[\Big\{n^{-1/2}\sum_{i=1}^n\big({\bf Y}_{1i}\big)_{\ell_1}\big({\bf Y}_{2i}\big)_{\ell_2}\Big\}^2\Big]={\rm E}\big[\big({\bf Y}_{11}\big)_{\ell_1}^2\big] \cdot {\rm E}\big[\big({\bf Y}_{21}\big)_{\ell_2}^2\big].
\end{align*}

Turn to the right-hand quantity in~\eqref{two}.  Since $\vspace{.7mm}\big(\big({\bf Y}\n_{11}\big)_{\ell_1},\ldots,\big({\bf Y}\n_{1n}\big)_{\ell_1}\big)$ and $\big(\big({\bf Y}\n_{21}\big)_{\ell_2},\ldots,\big({\bf Y}\n_{2n}\big)_{\ell_2}\big)$, for~$\ell_1=1,\dots,d_1$, $\ell_2=1,\dots,d_2$,  are independent,  
 we have, in view of Proposition~\ref{H2018}(ii) in Section~\ref{Propsec} and Theorem~2 in \citet{MR0044058}, 
\begin{align*}
&{\rm E}\Big[\Big\{n^{-1/2}\sum_{i=1}^n\big({\bf Y}\n_{1i}\big)_{\ell_1}\big({\bf Y}\n_{2i}\big)_{\ell_2}\Big\}^2\Big]
 =\frac{1}{n(n-1)}
\sum_{{\bf u}_1\in \mathfrak{G}_n^{d_1}} \Big(\big({\bf J}_1({\bf u}_1)\big)_{\ell_1}\Big)^2
\sum_{{\bf u}_2\in \mathfrak{G}_n^{d_2}} \Big(\big({\bf J}_2({\bf u}_2)\big)_{\ell_2}\Big)^2,
\end{align*}
where the right-hand side, by the properties of the grids $\mathfrak{G}^{d_k}_n$ and the fact that the score functions $J_k$, $k=1,2$, 
are continuous, square-integrable, and satisfy \eqref{scorevariance},
tends to
\begin{align*}
{\rm E}\big[\big({\bf J}_1({\bf V}_1)\big)_{\ell_1}^2\big] \cdot 
{\rm E}\big[\big({\bf J}_2({\bf V}_2)\big)_{\ell_2}^2\big]
={\rm E}\big[\big({\bf Y}_{11}\big)_{\ell_1}^2\big] \cdot 
 {\rm E}\big[\big({\bf Y}_{21}\big)_{\ell_2}^2\big]
\end{align*}
with ${\bf V}_k\sim {\rm U}_{d_k}$, $k=1,2$. 

Next, we obtain for \eqref{three} 
\begin{align}
&{\rm E}\Big[n^{-1/2}\sum_{i=1}^n\big({\bf Y}\n_{1i}\big)_{\ell_1}
                                 \big({\bf Y}\n_{2i}\big)_{\ell_2}\times
n^{-1/2}\sum_{i=1}^n\big({\bf Y}_{1i}\big)_{\ell_1}
                    \big({\bf Y}_{2i}\big)_{\ell_2}\Big]\nonumber\\
&=n^{-1}\Big\{\sum_{i=1}^n
\Big({\rm E}\big[\big({\bf Y}\n_{1i}\big)_{\ell_1}\big({\bf Y}_{1i}\big)_{\ell_1}\big]\Big)
\Big({\rm E}\big[\big({\bf Y}\n_{2i}\big)_{\ell_2}\big({\bf Y}_{2i}\big)_{\ell_2}\big]\Big)\nonumber\\
&\qquad+\sum_{i\ne i\pr}
\Big({\rm E}\big[\big({\bf Y}\n_{1i\pr}\big)_{\ell_1}\big({\bf Y}_{1i}\big)_{\ell_1}\big]\Big)
\Big({\rm E}\big[\big({\bf Y}\n_{2i\pr}\big)_{\ell_2}\big({\bf Y}_{2i}\big)_{\ell_2}\big]\Big)
\Big\}.\label{cross}
\end{align}
Since, for $k=1,2$, $\ell_k=1,\dots,d_k$,
\begin{align*}
{\rm E}\big[\big({\bf Y}\n_{ki}\big)_{\ell_k}\big({\bf Y}_{ki}\big)_{\ell_k}\big]
+\sum_{i\pr: i\pr\neq i}{\rm E}\big[\big({\bf Y}\n_{ki\pr}\big)_{\ell_k}\big({\bf Y}_{ki}\big)_{\ell_k}\big]
 &={\rm E}\Big[\big({\bf Y}_{ki}\big)_{\ell_k}\sum_{i\pr=1}^n\big({\bf Y}\n_{ki\pr}\big)_{\ell_k}\Big]\\
={\rm E}\Big[\big({\bf Y}_{ki}\big)_{\ell_k}\sum_{{\bf u}\in \mathfrak{G}_n^{d_k}}\big({\bf J}_k({\bf u})\big)_{\ell_k}\Big]
&={\rm E}\Big[\big({\bf Y}_{ki}\big)_{\ell_k}\Big]\Big(\sum_{{\bf u}\in \mathfrak{G}_n^{d_k}}{\bf J}_k({\bf u})\Big)_{\ell_k}
=0,
\end{align*}
and since, by exchangeability, for all $i\neq i\pr$, 
\begin{align*}
 {\rm E}\big[\big({\bf Y}\n_{ki}\big)_{\ell_k}\big({\bf Y}_{ki}\big)_{\ell_k}\big]
={\rm E}\big[\big({\bf Y}\n_{k1}\big)_{\ell_k}\big({\bf Y}_{k1}\big)_{\ell_k}\big]~~~\text{and}~~~
 {\rm E}\big[\big({\bf Y}\n_{ki\pr}\big)_{\ell_k}\big({\bf Y}_{ki}\big)_{\ell_k}\big]
={\rm E}\big[\big({\bf Y}\n_{k2}\big)_{\ell_k}\big({\bf Y}_{k1}\big)_{\ell_k}\big],
\end{align*}
we deduce that the right-hand side of \eqref{cross} equals
\begin{align}
&n^{-1}\Big\{n
\Big({\rm E}\big[\big({\bf Y}\n_{11}\big)_{\ell_1}\big({\bf Y}_{11}\big)_{\ell_1}\big]\Big)
\Big({\rm E}\big[\big({\bf Y}\n_{21}\big)_{\ell_2}\big({\bf Y}_{21}\big)_{\ell_2}\big]\Big)\nonumber\\
&\qquad\qquad+n(n-1)\Big(-(n-1)^{-1}{\rm E}
\big[\big({\bf Y}\n_{11}\big)_{\ell_1}\big({\bf Y}_{11}\big)_{\ell_1}\big]\Big)
 \Big(-(n-1)^{-1}{\rm E}
\big[\big({\bf Y}\n_{21}\big)_{\ell_2}\big({\bf Y}_{21}\big)_{\ell_2}\big]\Big)\Big\}\nonumber\\
&\qquad\qquad\qquad\qquad\qquad\qquad=\frac{n}{n-1}
{\rm E}\big[\big({\bf Y}\n_{11}\big)_{\ell_1}\big({\bf Y}_{11}\big)_{\ell_1}\big]\cdot
{\rm E}\big[\big({\bf Y}\n_{21}\big)_{\ell_2}\big({\bf Y}_{21}\big)_{\ell_2}\big].\label{cross2}
\end{align}
The Proof of Theorem 3.1 in \citet{del2018smooth} and the Proof of Proposition 3.3 in \citet{MR4255122} entails 
${\bf Y}\n_{k1}-{\bf Y}_{k1}\longrightarrow 0$ a.s.,  
while by \eqref{scorevariance}
\begin{align*}
\lim_{n\to\infty} {\rm E}\big[\Vert {\bf Y}\n_{k1}\Vert^2\big]
& 
={\rm E}\big[\Vert {\bf Y}_{k1}\Vert^2\big]
=\int_0^1 J_k^2(u)\rd u 
<\infty .
\end{align*}
It then follows (see, e.g., part {\it (iv)} of Theorem~5.7 in \citet[Chap.~3]{MR3701383})
that
$${\rm E}\big[\Vert {\bf Y}\n_{k1}-{\bf Y}_{k1}\Vert^2\big] \to 0,\quad k=1,2\ \text{ as $n\to\infty$}.$$
In particular,
 ${\rm E}\big[\big({\bf Y}\n_{k1}-{\bf Y}_{k1}\big)_{\ell_k}^2\big] \to 0$ 
and thus 
\begin{align*}
&{\rm E}\big[\big({\bf Y}\n_{k1}\big)_{\ell_k}\big({\bf Y}_{k1}\big)_{\ell_k}\big]
 =2^{-1}\Big({\rm E}\big[\big({\bf Y}\n_{k1}\big)_{\ell_k}^2\big]
+{\rm E}\big[\big({\bf Y}_{k1}\big)_{\ell_k}^2\big]
-{\rm E}\big[\big({\bf Y}\n_{k1}-{\bf Y}_{k1}\big)_{\ell_k}^2\big]\Big)
 \to {\rm E}\big[\big({\bf Y}_{k1}\big)_{\ell_k}^2\big].
\end{align*}
It follows that the right-hand side of \eqref{cross2} tends to~${\rm E}\big[\big({\bf Y}_{11}\big)_{\ell_1}^2\big] \cdot {\rm E}\big[\big({\bf Y}_{21}\big)_{\ell_2}^2\big]$, which concludes the proof   of the lemma for {\rm vec}$({\tenq{\mathbf W}}\,\!_{J}\n)$.

\medskip

\noindent  (b) The case of the Kendall matrix $\text{\rm vec}({\tenq{\mathbf W}}\,\!_{\text{\tiny\rm Kendall}}\n)$ is slightly different, although the arguments in the proof are quite similar. We consider the H\' ajek projection of U-statistics (see, e.g., Proof of Theorem~7.1 in \citet{MR0026294}) for
\begin{align*}
\Big( {{\mathbf W}}\,\!_{\text{\tiny\rm Kendall}}\n \Big)_{\ell_1 \ell_2}= {n \choose 2}^{-1} \sum_{i<i\pr}
\;&   \text{sign} \Big(\big({\bf F}_{1;\pms}({\bf X}_{1i}) - {\bf F}_{1;\pms}({\bf X}_{1i\pr}) \big)_{\ell_1}\Big)
\times\text{sign} \Big(\big({\bf F}_{2;\pms}({\bf X}_{2i}) - {\bf F}_{2;\pms}({\bf X}_{2i\pr}) \big)_{\ell_2}\Big).
\end{align*}
It follows from Application 9(d) in \citet{MR0026294} that
\begin{align}
\Big( {{\mathbf W}}\,\!_{\text{\tiny\rm Kendall}}\n \Big)_{\ell_1 \ell_2}
=\frac{2}{n}\sum_{i=1}^{n}
\Big\{ 2F_{1\ell_1}\Big(\big({\bf F}_{1;\pms}({\bf X}_{1i})\big)_{\ell_1} \Big) -1 \Big\}
\Big\{ 2F_{2\ell_2}\Big(\big({\bf F}_{2;\pms}({\bf X}_{2i})\big)_{\ell_2} \Big) -1 \Big\}
+o_{\text{\rm q.m.}}(n^{-1/2}),\label{HU}
\end{align}
where $F_{k\ell}$ denotes the cumulative distribution functions of $\big({\bf F}_{k;\pms}({\bf X}_{k}) \big)_{\ell}$, $k=1,2$, $\ell=1,\dots,d_k$.

We also have the H\' ajek projection of combinatorial statistics (see, e.g., page 242 of \citet{MR857081} and Chapter II.3.1 of \citet{MR0229351}) for
\begin{align*}
\Big( {\tenq{\mathbf W}}\,\!_{\text{\tiny\rm Kendall}}\n \Big)_{\ell_1 \ell_2}= {n \choose 2}^{-1} \sum_{i<i\pr}
\text{sign}\Big(\big({\bf F}\n_{1;\pms}({\bf X}_{1i}) - {\bf F}\n_{1;\pms}({\bf X}_{1i\pr}) \big)_{\ell_1}\Big)\times
\text{sign}\Big(\big({\bf F}\n_{2;\pms}({\bf X}_{2i}) - {\bf F}\n_{2;\pms}({\bf X}_{2i\pr}) \big)_{\ell_2}\Big),
\end{align*}
which implies
\begin{align}
\Big( {\tenq{\mathbf W}}\,\!_{\text{\tiny\rm Kendall}}\n \Big)_{\ell_1 \ell_2}
=\frac{2}{n}\sum_{i=1}^{n}
\Big\{ 2F\n_{1\ell_1;\text{\tiny\rm mid}}\Big(\big({\bf F}\n_{1;\pms}({\bf X}_{1i})\big)_{\ell_1} \Big) -1 \Big\}
\Big\{ 2F\n_{2\ell_2;\text{\tiny\rm mid}}\Big(\big({\bf F}\n_{2;\pms}({\bf X}_{2i})\big)_{\ell_2} \Big) -1 \Big\}
+o_{\text{\rm q.m.}}(n^{-1/2}),\label{HC}
\end{align}
where $F\n_{k\ell;\text{\tiny\rm mid}}$ denotes the {\it mid-cumulative distribution function}\footnote{The {\it mid-cumulative distribution function} of a random variable~$X$ is defined as $F_{\text{\tiny\rm mid}}(x)\coloneqq [{\rm P}(X\le x)+{\rm P}(X< x)]/2.$ } \citep{MR2185587} of $\big({\bf F}\n_{k;\pms}({\bf X}_{k}) \big)_{\ell}$, \linebreak $k=1,2$, $\ell=1,\dots,d_k$.

Finally, it follows along the same lines as in the proof for part~(1) that the difference between the right-hand sides of \eqref{HC} and \eqref{HU}  is~$o_{\text{\rm q.m.}}(n^{-1/2})$ as $n\to\infty$. This concludes that
\[
\Big({\tenq{\mathbf W}}\,\!_{\text{\tiny\rm Kendall}}\n\Big)_{\ell_1 \ell_2} -
\Big({{\mathbf W}}\,\!_{\text{\tiny\rm Kendall}}\n\Big)_{\ell_1 \ell_2}
=o_{\text{\rm q.m.}}(n^{-1/2}),
\]
which completes the proof.
\cqfd

\vspace{5mm}

\noindent {\bf Proof of Lemma~\ref{lem:lan}.} Write $\rP^{\bf X}_\delta$ for the distribution $\rP^{\bf X}_\delta(\rP_1,\rP_2,{\bf M}_1,{\bf M}_2)$ of $\bf X$. It follows from the quadratic mean differentiability of $f_1^{1/2}$ and~$f_2^{1/2}$ and the differentiability with respect to $\delta$ of   ${\bf M}_\delta
$ that, denoting by $[{\bf V}]_1$ and~$[{\bf V}]_2$, respectively,  the first~$d_1$ and last $d_2$ components of a $d$-dimensional vector~$\bf V$, 
\begin{equation*}
\delta\mapsto \Big( \frac{\rd \rP^{\bf X}_\delta}{\rd \mu_d}({\bf x})\Big)^{1/2} 
= \Big( \big\vert \text{det}\big({\bf M}_\delta\big)
\big\vert^{-1}\Big)^{1/2}
 f_1([{\bf M}^{-1}_{\delta}{\bf x}]_1)
f_2([{\bf M}^{-1}_{\delta}{\bf x}]_2)\quad \delta\in{\mathbb R},\  {\bf x}\in {\mathbb R}^{d}
\end{equation*}
also is differentiable in quadratic mean. The quadratic expansion of the log-likelihood ratio 
$$\delta\mapsto \log\frac{\rd\rP^{(n)}_{\delta + n^{-1/2}\tau}}{\rd\rP\n_\delta}({\bf X}\n)$$ follows (see, e.g., Theorem~12.2.3~(i) in \citet{MR2135927}), yielding, at $\delta=0$, the second-order asymptotic representation~\eqref{LANKon}.  The explicit forms of the central sequence~$\Delta\n({\bf X}\n)$ and the Fisher information~$\gamma^2$ for $\delta =0$ are obtained via elementary differentiation.  The asymptotic normality result for $\Delta\n({\bf X}\n)$ follows from part (ii) of the same Theorem~12.2.3. 
 \cqfd 

\vspace{5mm}

\noindent {\bf Proof of Theorem~\ref{thm:J}.} 
We only give the proof for ${\tenq T}\,\!_{J}\n$;
 the proof for ${\tenq T}\,\!_{\text{\tiny\rm Kendall}}\n$ is similar and hence is omitted. 
Applying the multivariate central limit theorem \citep[Equation~(18.24)]{MR855460}
to the asymptotic form of~$\Lambda\n({\bf X}\n)$ (see Lemma~3.1), 
we obtain, under the null hypo\-thesis~($\delta = 0$),
\begin{equation*}
\Big(n^{1/2}{\rm vec}({{\mathbf W}}\,\!_{J}\n), \Lambda\n({\bf X}\n)\Big)
\stackrel{}{\rightsquigarrow}
N_{d_1d_2+1}\Bigg(
\Bigg(\begin{matrix}
{\bf 0}_{d_1d_2}\\
-\frac{1}{2}\tau^2\gamma^2
\end{matrix}\Bigg),
\Bigg(\begin{matrix}
\sigma^2_J{\bf I}_{d_1d_2} & \tau{\bf v}\\
\tau{\bf v}^{\prime} & \tau^2\gamma^2
\end{matrix}\Bigg)
\Bigg)\qquad\text{as $n\to\infty$}
\end{equation*}
where $\sigma^2_J\coloneqq {\sigma^{2}_{J_1}\sigma^{2}_{J_2}}/({d_1d_2})$ and 
\begin{align*}
{\bf v}&\coloneqq {\rm Cov}_{H_0}\Big[{\rm vec}\Big(
{\bf J}_1({\bf F}_{1;\pms}({\bf X}_{1}))
{\bf J}_2({\bf F}_{2;\pms}({\bf X}_{2}))\pr\Big),\hspace{30mm}\\
&\hspace{30mm}{\bf X}\pr_{1}{\bf M}_2\pr{\boldsymbol\varphi}_2({\bf X}_{2}) + 
{\bf X}\pr_{2}{\bf M}_1\pr{\boldsymbol\varphi}_1({\bf X}_{1}) 
 -\Big({\bf X}\pr_{1}{\boldsymbol\varphi}_1({\bf X}_{1}) -d_1
\Big)
-\Big({\bf X}\pr_{2}{\boldsymbol\varphi}_2({\bf X}_{2}) -d_2
\Big)
\Big]\\
&\;={\rm E}_{H_0}\Big[{\rm vec}\Big(
{\bf J}_1({\bf F}_{1;\pms}({\bf X}_{1}))
{\bf J}_2({\bf F}_{2;\pms}({\bf X}_{2}))\pr\Big)\Big(
{\bf X}\pr_{1}{\bf M}_2\pr{\boldsymbol\varphi}_2({\bf X}_{2}) + 
{\bf X}\pr_{2}{\bf M}_1\pr{\boldsymbol\varphi}_1({\bf X}_{1})\Big)\Big] \\
&\;={\rm vec}\Big({\rm E}_{H_0}\Big[
{\bf J}_1({\bf F}_{1;\pms}({\bf X}_{1})){\bf R}
{\bf J}_2({\bf F}_{2;\pms}({\bf X}_{2}))\pr\Big]\Big).
\end{align*}
Thus, by Lemma~\ref{lem:lan},
\begin{equation*}
\Big(n^{1/2}{\rm vec}({\tenq{\mathbf W}}\,\!_{J}\n), \Lambda\n({\bf X}\n)\Big)
\stackrel{}{\rightsquigarrow}
N_{d_1d_2+1}\Bigg(
\Bigg(\begin{matrix}
{\bf 0}_{d_1d_2}\\
-\frac{1}{2}\tau^2\gamma^2
\end{matrix}\Bigg),
\Bigg(\begin{matrix}
\sigma^2_J{\bf I}_{d_1d_2} & \tau{\bf v}\\
\tau{\bf v}^{\prime} & \tau^2\gamma^2
\end{matrix}\Bigg)
\Bigg)\qquad\text{as $n\to\infty$}.
\end{equation*}
Le~Cam's third lemma 
\citep[Chapter~VI.1.4]{MR0229351} then yields, under local alterna-\linebreak tives~($\delta = n^{-1/2}\tau$),
\[n^{1/2}{\rm vec}({\tenq{\mathbf W}}\,\!_{J}\n)
\stackrel{}{\rightsquigarrow}
N_{d_1d_2}(\tau{\bf v},\sigma^2_J{\bf I}_{d_1d_2})\qquad\text{as $n\to\infty$}. 
\]
The result follows.  \cqfd 

\vspace{5mm}

\noindent {\bf Proof of Proposition~\ref{prop:Pitman}.} 
Denoting ${\bf Y}_{k}:={\boldsymbol{\Sigma}}_{k}^{-1/2}{\bf X}_{k}$, 
and recalling $\widetilde F_k$ as the cumulative distribution function of $\Vert {\bf Y}_{k}\Vert$, $k=1,2$,
direct computation yields
\begin{align*}
&{\rm E}_{H_0}\Big[{\bf J}_1({\bf F}_{1;\pms} ({\bf X}_1))
\Big({\bf X}_{1}\pr{\bf M}_2\pr{\boldsymbol\varphi}_2({\bf X}_{2})
\Big) {\bf J}_2({\bf F}_{2;\pms} ({\bf X}_2))\pr\Big] \\
&\quad ={\rm E}_{H_0}\Big[{\bf J}_1\Big(\frac{{\bf Y}_{1}}{\Vert {\bf Y}_{1}\Vert}\widetilde F_1(\Vert {\bf Y}_{1}\Vert)\Big)
\Big({\bf Y}_{1}\pr{\boldsymbol{\Sigma}}_{1}^{1/2}
 {\bf M}_2\pr{\boldsymbol\varphi}_2({\boldsymbol{\Sigma}}_{2}^{1/2}{\bf Y}_{2})
\Big) {\bf J}_2\Big(\frac{{\bf Y}_{2}}{\Vert {\bf Y}_{2}\Vert}\widetilde F_2(\Vert {\bf Y}_{2}\Vert)\Big)\pr\Big] \\
&\quad ={\rm E}_{H_0}\Big[\frac{{\bf Y}_{1}}{\Vert {\bf Y}_{1}\Vert}J_1\Big(\widetilde F_1(\Vert {\bf Y}_{1}\Vert)\Big)
\Big({\bf Y}_{1}\pr{\boldsymbol{\Sigma}}_{1}^{1/2}
 {\bf M}_2\pr{\boldsymbol{\Sigma}}_{2}^{-1/2}\frac{{\bf Y}_2}{\Vert {\bf Y}_{2}\Vert}\rho_2(\Vert {\bf Y}_{2}\Vert)
\Big) J_2\big(\widetilde F_2(\Vert {\bf Y}_{2}\Vert)\big)\frac{{\bf Y}_{2}\pr}{\Vert {\bf Y}_{2}\Vert}\Big]\\
&\quad ={\rm E}_{H_0}\Big[\frac{{\bf Y}_{1}}{\Vert {\bf Y}_{1}\Vert}J_1\Big(\widetilde F_1(\Vert {\bf Y}_{1}\Vert)\Big){\bf Y}_{1}\pr\Big]
{\boldsymbol{\Sigma}}_{1}^{1/2}{\bf M}_2\pr{\boldsymbol{\Sigma}}_{2}^{-1/2}
{\rm E}_{H_0}\Big[\frac{{\bf Y}_2}{\Vert {\bf Y}_{2}\Vert}\rho_2(\Vert {\bf Y}_{2}\Vert)
 J_2\big(\widetilde F_2(\Vert {\bf Y}_{2}\Vert)\big)\frac{{\bf Y}_{2}\pr}{\Vert {\bf Y}_{2}\Vert}\Big].
\end{align*}
Write $r_k:=\Vert {\bf Y}_{k}\Vert$ and ${\bf U}_{k}:={\bf Y}_{k}\big/\Vert {\bf Y}_{k}\Vert$, $k=1,2$.  
We have, by the independence between $r_k$ and ${\bf U}_{k}$ \citep[Theorem~2.3]{MR1071174},
\begin{align*}
{\rm E}_{H_0}\Big[\frac{{\bf Y}_{1}}{\Vert {\bf Y}_{1}\Vert}J_1\Big(\widetilde F_1(\Vert {\bf Y}_{1}\Vert)\Big){\bf Y}_{1}\pr\Big]
&={\rm E}_{H_0}\Big[{\bf U}_{1} J_1\Big(\widetilde F_1(r_1)\Big)r_1{\bf U}_{1}\pr\Big]\\
&={\rm E}_{H_0}\Big[J_1\Big(\widetilde F_1(r_1)\Big)r_1\Big] \cdot {\rm E}_{H_0}\Big[{\bf U}_{1} {\bf U}_{1}\pr\Big]
=D_1 d_1^{-1}{\bf I}_{d_1},
\end{align*}
and similarly,
\[
{\rm E}_{H_0}\Big[\frac{{\bf Y}_2}{\Vert {\bf Y}_{2}\Vert}\rho_2(\Vert {\bf Y}_{2}\Vert)
 J_2\Big(\widetilde F_2(\Vert {\bf Y}_{2}\Vert)\Big)\frac{{\bf Y}_{2}\pr}{\Vert {\bf Y}_{2}\Vert}\Big]
={\rm E}_{H_0}\Big[\rho_2(r_2)
 J_2\Big(\widetilde F_2(r_2)\Big)\Big] \cdot {\rm E}_{H_0}\Big[{\bf U}_2{\bf U}_{2}\pr\Big]
=C_2 d_2^{-1}{\bf I}_{d_2},
\]
where ${\rm E}_{H_0}\Big[{\bf U}_{k}{\bf U}_{k}\pr\Big]=d_k^{-1}{\bf I}_{d_k}$, $k=1,2$; see Theorem~2.7 in \citet{MR1071174}. Therefore,
\begin{equation}\label{add1}
{\rm E}_{H_0}\Big[{\bf J}_1({\bf F}_{1;\pms} ({\bf X}_1))
\Big({\bf X}_{1}\pr{\bf M}_2\pr{\boldsymbol\varphi}_2({\bf X}_{2})
\Big) {\bf J}_2({\bf F}_{2;\pms} ({\bf X}_2))\pr\Big] =d_1^{-1}d_2^{-1}D_1C_2{\boldsymbol{\Sigma}}_{1}^{1/2}{\bf M}_2\pr{\boldsymbol{\Sigma}}_2^{-1/2},
\end{equation}
and similarly,
\begin{equation}\label{add2}
{\rm E}_{H_0}\Big[{\bf J}_1({\bf F}_{1;\pms} ({\bf X}_1))
\Big({\bf X}_{2}\pr{\bf M}_1\pr{\boldsymbol\varphi}_1({\bf X}_{1})
\Big) {\bf J}_2({\bf F}_{2;\pms} ({\bf X}_2))\pr\Big] =d_1^{-1}d_2^{-1}D_2C_1{\boldsymbol{\Sigma}}_{1}^{-1/2}{\bf M}_1{\boldsymbol{\Sigma}}_2^{1/2}.
\end{equation}
Plugging \eqref{add1} and \eqref{add2} into \eqref{still} yields that, under the alterna\-tive ${\bf X}\sim {\rm P}_{n^{-1/2}\tau}^{\,\text{\rm{\tiny ell}}}(\phi_1,\phi_2,{\boldsymbol\Sigma}_1, {\boldsymbol\Sigma}_2; {\bf M}_1, {\bf M}_2)$, the limiting distribution of the test statistic ${\tenq T}\,\!_{J}\n$ is noncentral chi-square with $d_1d_2$ degrees of freedom and noncentrality parameter 
 \begin{equation*}
\frac{\tau^2 }{d_1d_2\sigma^{2}_{J_1}\sigma^{2}_{J_2}}
\Big\Vert  D_1C_2{\boldsymbol{\Sigma}}_{1}^{1/2}{\bf M}_2\pr{\boldsymbol{\Sigma}}_2^{-1/2}
+ D_2C_1{\boldsymbol{\Sigma}}_{1}^{-1/2}{\bf M}_1{\boldsymbol{\Sigma}}_2^{1/2} \Big\Vert^2_{\mathrm F},
 \end{equation*}
and recall that the limiting distribution of Wilks' statistic is noncentral chi-square with~$d_1d_2$ degrees of freedom and noncentrality parameter of \eqref{still2}.  
The first result~\eqref{AREell} then follows from the fact that  (see, e.g., \citet[Equation~(5)]{MR83228}) when two test statistics are asymptotically noncentral chi-squared distributed under a local sequence of alternatives (here,~$\rP\n _{n^{-1/2}\tau}$),  their Pitman asymptotic relative efficiencies are obtained as the ratios of their noncentrality parameters.  
Claims (i) and (ii)  follow from the proofs of Propositions 1 and~2 in \citet{MR2462206}; see also Theorem 1 in \citet{MR2160622} and Proposition 7 in \citet{MR1963662}. 
  \cqfd

\begin{acks}[Acknowledgments]
The authors would like to thank two anonymous referees, an anonymous Associate
Editor and the Editor Davy Paindaveine for their constructive comments that improved the
quality of this paper.
\end{acks}

\begin{funding}
Hongjian Shi and Mathias Drton were supported by the European Research Council (ERC) under the European Union's Horizon 2020 research and innovation programme (grant agreement No 883818). Marc Hallin acknowledges the support of the Czech Science Foundation grant GA\v{C}R22036365. Fang~Han was supported by the United States NSF grants DMS-1712536 and SES-2019363.
\end{funding}

\begin{supplement}
\stitle{Supplement to ``Distribution-free tests of multivariate independence based on center-outward quadrant, Spearman, Kendall, and van der Waerden statistics''}
\sdescription{The supplement consists of auxiliary results and numerical experiments.}
\end{supplement}


\bibliographystyle{imsart-nameyear} 
\bibliography{AMS}       



\section*{Supplementary material}

**Proof of Proposition 4.4.** Writing $\mathrm{E}_{H_0}$ for expectations under the null ($\delta = 0$), direct computation yields

$$\mathrm{E}_{H_0}\Big[\mathbf{J}_1(\mathbf{F}_{1;\pm}(\mathbf{X}_1))\Big(\mathbf{X}_1'\mathbf{M}_2'\boldsymbol{\varphi}_2(\mathbf{X}_2)\Big)\mathbf{J}_2(\mathbf{F}_{2;\pm}(\mathbf{X}_2))'\Big]$$
$$= \mathrm{E}_{H_0}\Big[\mathbf{J}_1\Big(\frac{\mathbf{Y}_1}{\|\mathbf{Y}_1\|}\tilde{F}_1(\|\mathbf{Y}_1\|)\Big)\Big(\mathbf{Y}_1'\boldsymbol{\Sigma}_1^{1/2}\mathbf{M}_2'\boldsymbol{\varphi}_2(\boldsymbol{\Sigma}_2^{1/2}\mathbf{Y}_2)\Big)\mathbf{J}_2\Big(\frac{\mathbf{Y}_2}{\|\mathbf{Y}_2\|}\tilde{F}_2(\|\mathbf{Y}_2\|)\Big)'\Big]$$
$$= \mathrm{E}_{H_0}\Big[\frac{\mathbf{Y}_1}{\|\mathbf{Y}_1\|}J_1\Big(\tilde{F}_1(\|\mathbf{Y}_1\|)\Big)\Big(\mathbf{Y}_1'\boldsymbol{\Sigma}_1^{1/2}\mathbf{M}_2'\boldsymbol{\varphi}_2(\boldsymbol{\Sigma}_2^{1/2}\mathbf{Y}_2)\Big)J_2\Big(\tilde{F}_2(\|\mathbf{Y}_2\|)\Big)\frac{\mathbf{Y}_2'}{\|\mathbf{Y}_2\|}\Big]$$
$$= D_1 C_2 \boldsymbol{\Sigma}_1^{1/2} \mathbf{M}_2' \boldsymbol{\Sigma}_2^{-1/2}.$$

The first result (16) then follows from the fact that (see, e.g., Hannan (1956, Equation (5))) when two test statistics are asymptotically noncentral chi-squared distributed under a local sequence of alternatives (here, $\mathrm{P}_{n^{-1/2}\tau}^{(n)}$), their Pitman asymptotic relative efficiencies are obtained as the ratios of their noncentrality parameters. Claims (i) and (ii) follow from the proofs of Propositions 1 and 2 in Hallin and Paindaveine (2008); see also Theorem 1 in Paindaveine (2004) and Proposition 7 in Hallin and Paindaveine (2002). □

## S.2 Further numerical results

### S.2.1 Further analysis of the motivating example

Supporting our claim that forty observations three months apart can be considered as independent, Figure S.1 provides correlogram and cross-correlogram plots for the series in the motivating example of Section 1.3. Observe that, except for lag zero, essentially all auto-and cross-correlations are non-significant at level 5%, and that no clear seasonality pattern emerges. Therefore, treating these observations as independent and identically distributed is unlikely to affect the validity of the tests.

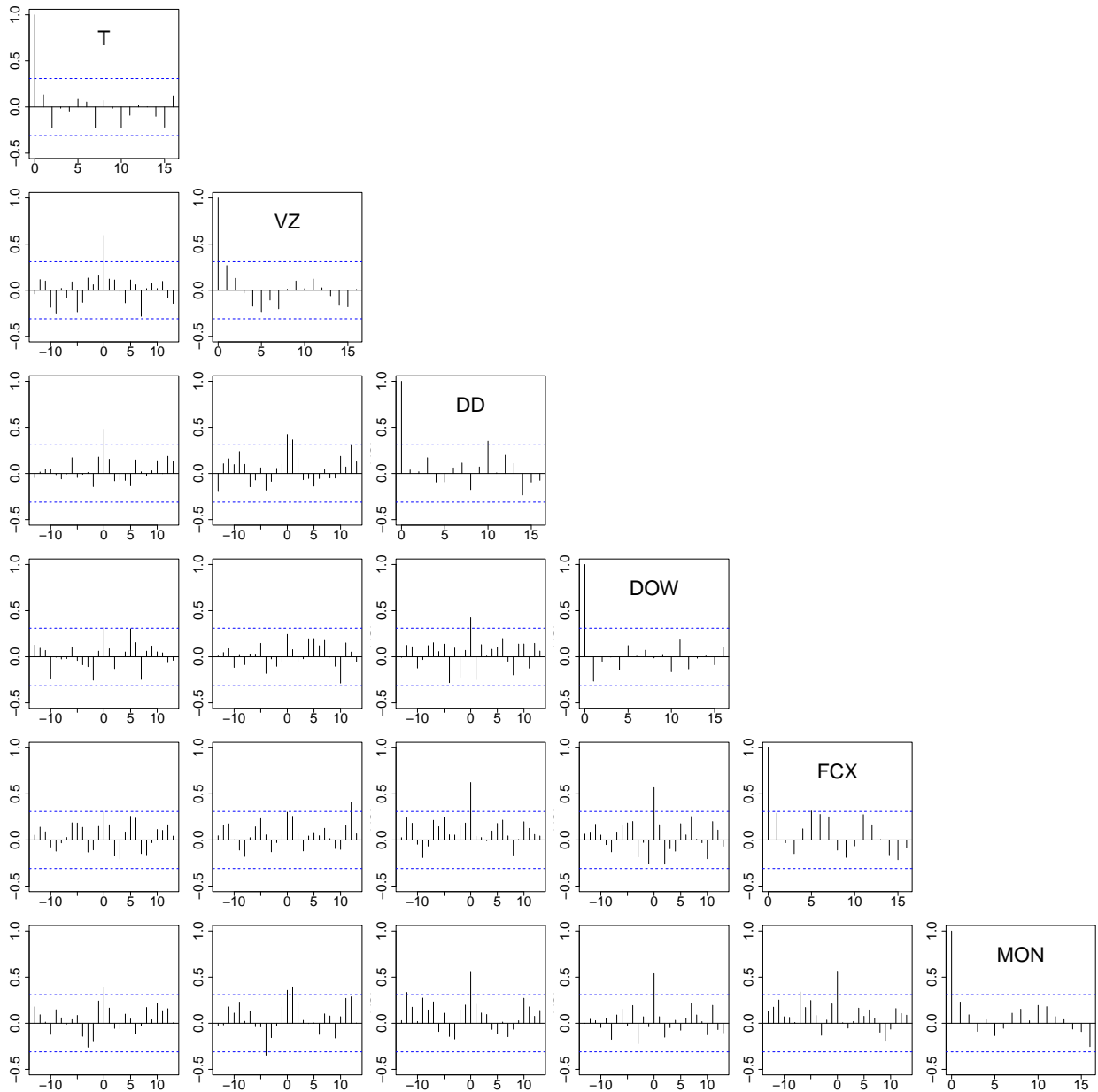

Figure S.1: Autocorrelation (on the diagonal) and cross-correlation plots (below the diagonal), with 95% confidence bands, for the forty observations of daily log returns of the stocks "AT&T Inc [T]", "Verizon Communications [VZ]," "Du Pont (E.I.) [DD]," "Dow Chemical [DOW]," "Freeport-McMoran Cp & Gld [FCX]," and "Monsanto Co. [MON]."

4

### S.2.2 Finite-sample size of Wilks' likelihood ratio tests

Likelihood ratio tests, although asymptotically valid, often suffer from severe size inflation when the observations are not Gaussian, and Wilks' test is no exception. While this is hardly surprising, the magnitude of this inflation is typically underestimated; Table S.1 provides a brief but quite alarming simulation-based illustration of the phenomenon for various dimensions $d_1$ and $d_2$. Thanks to distribution-freeness, our rank-based tests are not affected, which constitutes another strong argument in their favor.

The simulated data are samples (size $n = 100$) of independent copies of random vectors of the form $\mathbf{X} = (\mathbf{X}_1', \mathbf{X}_2')'$ where $\mathbf{X}_1 = \exp(\mathbf{Z}_1)$ (values in $\mathbb{R}^{d_1}$), $\mathbf{X}_2 = \exp(\mathbf{Z}_2)$ (values in $\mathbb{R}^{d_2}$), with $(\mathbf{Z}_1, \mathbf{Z}_2) \sim N(\mathbf{0}, \sigma^2 \mathbf{I}_{d_1+d_2})$ and $\exp(\mathbf{Z}) := (\exp Z_1, \ldots, \exp(Z_d))'$. The empirical sizes of Wilks' test with nominal level 0.05, based on 10,000 simulations, are reported in Table S.1 for dimensions $d_1 = d_2 \in \{2, 3, 5, 7\}$ and $\sigma \in \{0.5, 1.0, 1.5, 2.0\}$. Note that $\sigma$, a standard error for the Gaussian $\mathbf{Z}_1$ and $\mathbf{Z}_2$, is a skewness parameter for the lognormal components of $\mathbf{X}_1$ and $\mathbf{X}_2$. As one might expect, the size spectacularly grows beyond the nominal level 0.05 for increasing skewness $\sigma$ and/or increasing dimensions $d_1$ and $d_2$, reaching 0.317 for $d_1 = d_2 = 7$ and $\sigma = 2.0$! Wilks' test, thus, is highly unreliable under non-Gaussian distributions and its conclusions should be interpreted with highest precaution.

Table S.1: Sizes of Wilks' LRT for log-normal distributions at nominal significance level 0.05 for various dimensions $d_1$ and $d_2$ and various values of the skewness parameter $\sigma$.

|  | $\sigma = 0.5$ | $\sigma = 1.0$ | $\sigma = 1.5$ | $\sigma = 2.0$ |
|---|---|---|---|---|
| $(d_1, d_2) = (2, 2)$ | 0.059 | 0.079 | 0.101 | 0.102 |
| $(d_1, d_2) = (3, 3)$ | 0.067 | 0.101 | 0.134 | 0.155 |
| $(d_1, d_2) = (5, 5)$ | 0.095 | 0.154 | 0.211 | 0.243 |
| $(d_1, d_2) = (7, 7)$ | 0.133 | 0.197 | 0.272 | 0.317 |

### S.2.3 The elliptical Hodges–Lehmann lower bound $\Omega(d_1, d_2)$

Table S.2 provides, for dimensions $d_1, d_2 \leq 10$, numerical values of the elliptical Hodges–Lehmann lower bounds $\Omega(d_1, d_2)$ (Proposition 4.4) for the AREs of Spearman center-outward rank tests (based on Wilcoxon scores) with respect to Wilks' test.

The bound, which, with a value of 0.913, is quite high for $d_1 = d_2 = 2$, slowly decreases as the dimensions increase. Recall that an ARE of 0.913 means that, in the worst cas, Wilks requires $0.913n$ observations to achieve the same performance as Spearman with $n$ observations. The worst efficiency loss, thus, is minimal; on the other hand, the same ARE has no upper bound, and Spearman, potentially, can be arbitrarily more powerful than Wilks.

### S.2.4 Additional simulation results

The data are generated as a sample of $n$ independent copies of the $(d_1 + d_2)$-dimensional random vector $\mathbf{X} = (\mathbf{X}_1', \mathbf{X}_2')'$ under the model described in Section 5. The following null distributions were considered with $d_1 = d_2$; note that (e), (h), and (i) are not elliptical, while the heavy-tailed (f)-(i)

Table S.2: Some numerical values of the Hodges–Lehmann lower bound $\Omega(d_1, d_2)$ for the Pitman asymptotic relative efficiency, against elliptical Konijn alternatives, of the center-outward Spearman (Wilcoxon score) test with respect to Wilks', for dimensions $d_1, d_2 \leq 10$.

| $d_1$ \ $d_2$ | 1 | 2 | 3 | 4 | 5 | 6 | 7 | 8 | 9 | 10 |
|---|---|---|---|---|---|---|---|---|---|---|
| 1 | 0.856 | 0.884 | 0.867 | 0.850 | 0.837 | 0.826 | 0.817 | 0.809 | 0.803 | 0.797 |
| 2 | 0.884 | 0.913 | 0.895 | 0.878 | 0.864 | 0.853 | 0.844 | 0.836 | 0.829 | 0.823 |
| 3 | 0.867 | 0.895 | 0.878 | 0.861 | 0.847 | 0.836 | 0.827 | 0.819 | 0.813 | 0.807 |
| 4 | 0.850 | 0.878 | 0.861 | 0.845 | 0.831 | 0.820 | 0.811 | 0.804 | 0.797 | 0.792 |
| 5 | 0.837 | 0.864 | 0.847 | 0.831 | 0.818 | 0.807 | 0.799 | 0.791 | 0.785 | 0.779 |
| 6 | 0.826 | 0.853 | 0.836 | 0.820 | 0.807 | 0.797 | 0.788 | 0.781 | 0.775 | 0.769 |
| 7 | 0.817 | 0.844 | 0.827 | 0.811 | 0.799 | 0.788 | 0.779 | 0.772 | 0.766 | 0.761 |
| 8 | 0.809 | 0.836 | 0.819 | 0.804 | 0.791 | 0.781 | 0.772 | 0.765 | 0.759 | 0.754 |
| 9 | 0.803 | 0.829 | 0.813 | 0.797 | 0.785 | 0.775 | 0.766 | 0.759 | 0.753 | 0.748 |
| 10 | 0.797 | 0.823 | 0.807 | 0.792 | 0.779 | 0.769 | 0.761 | 0.754 | 0.748 | 0.742 |

fail to satisfy the conditions for the validity of Wilks and other non-distribution-free tests, and the computation of AREs.

(e) (Gaussian mixture) $\mathbf{X}_1^*, \mathbf{X}_2^* \overset{\text{i.i.d.}}{\sim} 0.5 N(\mathbf{0}, \mathbf{I}_{d_0}) + 0.5 N(\mathbf{1}, \mathbf{I}_{d_0})$;

(f) (heavy-tailed elliptical) $\mathbf{X}_1^*, \mathbf{X}_2^* \overset{\text{i.i.d.}}{\sim} t_2(\mathbf{0}, \mathbf{I}_{d_0})$, where $t_\nu(\boldsymbol{\mu}, \boldsymbol{\Sigma})$ with $\boldsymbol{\mu} \in \mathbb{R}^p$, $\boldsymbol{\Sigma} \in \mathbb{R}^{p \times p}$ denotes the distribution with density

$$\frac{\Gamma[(\nu+p)/2]}{\Gamma(\nu/2)\nu^{p/2}\pi^{p/2} |\boldsymbol{\Sigma}|^{1/2}} \left[1 + \frac{1}{\nu}(\mathbf{x} - \boldsymbol{\mu})^T \boldsymbol{\Sigma}^{-1}(\mathbf{x} - \boldsymbol{\mu})\right]^{-(\nu+p)/2};$$

(g) (heavy-tailed elliptical Cauchy) $\mathbf{X}_1^*, \mathbf{X}_2^* \overset{\text{i.i.d.}}{\sim} t_1(\mathbf{0}, \mathbf{I}_{d_0})$;

(h) (heavy-tailed non-elliptical) $\mathbf{X}_{1i}^*, \mathbf{X}_{2i}^* \overset{\text{i.i.d.}}{\sim} t_2(0,1)$ for $i = 1, \ldots, d_0$;

(i) (heavy-tailed non-elliptical marginal Cauchy) $\mathbf{X}_{1i}^*, \mathbf{X}_{2i}^* \overset{\text{i.i.d.}}{\sim} t_1(0,1)$ for $i = 1, \ldots, d_0$.

Tables S.3–S.7 report empirical powers (rejection frequencies) of the seven tests considered in Section 5: (i) the center-outward sign (quadrant) test, (ii) the center-outward Spearman test, (iii) the center-outward Kendall test, (iv) the center-outward van der Waerden test, (v) the center-outward Wilcoxon score version of the distance covariance test ("dCov"), (vi) the center-outward van der Waerden score version of the distance covariance test ("dCvdW"), (vii) Wilks' likelihood ratio test, and (viii) Székely and Rizzo's permutation-based distance covariance test. These rejection frequencies are based on 1,000 simulations with nominal significance level $\alpha = 0.05$, dimensions $d_1 = d_2 \in \{2, 3, 5, 7\}$, sample size $n \in \{432, 864, 1728\}$, and parameter value $\delta^{(n)} = n^{-1/2}\tau$ with $\tau \in \{0, 0.2, 0.4, 0.6, 0.8\}$.

The results amply confirm the findings of Section 5. The center-outward van der Waerden test has the best power except when the true underlying distribution is Gaussian (case (a) of Section 5) or nearly Gaussian (case (e)); Wilks then, quite unsurprisingly, prevails. The performances of the center-outward Spearman and van der Waerden score tests, and the two distance covariance tests are comparable, while the center-outward sign and Kendall tests perform more poorly in some cases, which is not unexpected since the center-outward ranks and signs are not efficiently used.

Table S.3 shows the performance of the seven competing tests under normal mixtures, which are

non-elliptical but nearly Gaussian. Wilks' test has the best finite-sample performance compared to all the center-outward sign- and rank-based tests with, however, the center-outward van der Waerden test achieving similar performance for $n$ large enough.

Although cases (f) and (g) are not covered by Proposition 4.4 due to infinite second-order moments, Tables S.4 and S.5 indicate that the center-outward van der Waerden test remains best while, due to heavy-tailedness, Wilks's test badly suffers from severe over-rejection.

Tables S.6–S.7 deal with non-elliptical heavy-tailed distributions with infinite second-order moments. The center-outward van der Waerden test still performs best while Wilks suffers, as in the multivariate $t$ cases (f) and (g) (Tables S.4 and S.5), from severe size inflation due to heavy-tailedness. The only difference between case (g) (multivariate Cauchy; Table S.5) and case (i) (Cauchy marginals; Table S.7), is that the center-outward Kendall test performs poorly under the latter as the dimension increases.

To conclude this section, we explore mixture alternatives under which the density of $\mathbf{X} = (\mathbf{X}_1', \mathbf{X}_2')'$ is of the form
$$f_\delta := (1-\delta)f_1 f_2 + \delta f_0 \qquad \delta \in [0,1].$$
where $f_1$ and $f_2$ are the densities of the null distributions $P_1$ and $P_2$ of $\mathbf{X}_1$ and $\mathbf{X}_2$, respectively, and $f_0$ is the density of a distribution $P_0$ over $\mathbb{R}^{d_1+d_2}$ which does not factorize into the product of its marginals over $\mathbb{R}^{d_1}$ and $\mathbb{R}^{d_2}$. Clearly, the value $\delta = 0$ yields mutually independent vectors $\mathbf{X}_1$ and $\mathbf{X}_2$. For given $P_1$, $P_2$, and $P_0$, the distribution $P^\mathbf{X}$ of $\mathbf{X}$ belongs to a one-parameter family indexed by $\delta \in [0,1]$. LAN and ARE results for these families could be obtained along the same lines as in Section 4; details are left to the reader. Rather, we performed simulations similar to those in Sections 5 and S.2.4, with the following distributions $P_1$, $P_2$, and $P_0$ and the same values of $n$, $d_1 = d_2 = d_0$, and $\delta^{(n)} = n^{-1/2}\tau$:

(j) (standard Gaussian) $P_1 = P_2 = N(\mathbf{0}, \mathbf{I}_{d_0})$, and $P_0 = N(\mathbf{0}, \mathbf{A})$ where $\mathbf{A}$ is a symmetric $2d_0 \times 2d_0$ matrix with entries
$$A_{ij} = A_{ji} = \begin{cases} 1 & i = j, \\ 0.5 & i = j + d_0 \text{ or } j = i + d_0, \\ 0 & \text{otherwise;} \end{cases}$$

(k) (heavy-tailed elliptical) $P_1 = P_2 = t_3(\mathbf{0}, \mathbf{I}_{d_0})$, and $P_0 = t_3(\mathbf{0}, \mathbf{A})$ (where $t_3(\mathbf{0}, \mathbf{M})$ stands for the elliptical distribution with scatter matrix $\mathbf{M}$ and radial density $t_3$ centered at $\mathbf{0}$).

Tables S.8–S.9 report empirical powers (rejection frequencies) of the same seven tests as above, based on 1,000 simulations with nominal significance level $\alpha = 0.05$, sample size $n \in \{432, 864, 1728\}$, dimensions $d_1 = d_2 = d_0 \in \{2, 3, 5, 7\}$, and parameter values $\delta^{(n)} = n^{-1/2}\tau$, $\tau \in \{0, 1, 2, 3, 4\}$. When the sample size is sufficiently large, the center-outward van der Waerden test exhibits superior performance compared to, or is closely aligned with, both Wilks' and the distance covariance permutation-based test.

Table S.3: Empirical powers, under the Gaussian mixture (e), of *(i)* the center-outward sign (quadrant) test, *(ii)* the center-outward Spearman test, *(iii)* the center-outward Kendall test, *(iv)* the center-outward van der Waerden test, *(v)* the center-outward Wilcoxon score version of the distance covariance test ("dCov"), *(vi)* the center-outward van der Waerden score version of the distance covariance test ("dCvdW"), *(vii)* Wilks' likelihood ratio test, and *(viii)* Székely and Rizzo's permutation-based distance covariance test. Based on 1,000 replications of a sample with dimensions $d_1 = d_2 \in \{2, 3, 5, 7\}$, size $n \in \{432, 864, 1728\}$, and parameter values $\delta^{(n)} = n^{-1/2}\tau$, $\tau \in \{0, 0.2, 0.4, 0.6, 0.8\}$.

| | $n$ | 432 | | | | | 864 | | | | | 1728 | | | | |
|---|---|---|---|---|---|---|---|---|---|---|---|---|---|---|---|---|
| | $\tau$ | 0 | 0.2 | 0.4 | 0.6 | 0.8 | 0 | 0.2 | 0.4 | 0.6 | 0.8 | 0 | 0.2 | 0.4 | 0.6 | 0.8 |
| $(d_1, d_2) = (2,2)$ | sign | 0.050 | 0.052 | 0.093 | 0.145 | 0.247 | 0.049 | 0.058 | 0.095 | 0.161 | 0.262 | 0.045 | 0.063 | 0.093 | 0.157 | 0.255 |
| | Spearman | 0.045 | 0.068 | 0.109 | 0.217 | 0.389 | 0.055 | 0.069 | 0.129 | 0.230 | 0.403 | 0.047 | 0.061 | 0.114 | 0.235 | 0.386 |
| | Kendall | 0.049 | 0.057 | 0.111 | 0.213 | 0.379 | 0.051 | 0.065 | 0.115 | 0.219 | 0.380 | 0.045 | 0.063 | 0.116 | 0.225 | 0.355 |
| | vdW | 0.043 | 0.066 | 0.113 | 0.222 | 0.394 | 0.049 | 0.067 | 0.121 | 0.235 | 0.418 | 0.041 | 0.060 | 0.117 | 0.243 | 0.393 |
| | dCov | 0.050 | 0.059 | 0.098 | 0.201 | 0.367 | 0.047 | 0.067 | 0.120 | 0.224 | 0.373 | 0.049 | 0.061 | 0.103 | 0.218 | 0.345 |
| | dCvdW | 0.050 | 0.057 | 0.109 | 0.203 | 0.374 | 0.046 | 0.064 | 0.127 | 0.219 | 0.376 | 0.045 | 0.062 | 0.107 | 0.216 | 0.343 |
| | Wilks | 0.047 | 0.064 | 0.118 | 0.248 | 0.426 | 0.055 | 0.065 | 0.130 | 0.245 | 0.430 | 0.051 | 0.058 | 0.116 | 0.235 | 0.405 |
| | SzR | 0.055 | 0.065 | 0.121 | 0.226 | 0.356 | 0.030 | 0.071 | 0.111 | 0.214 | 0.368 | 0.051 | 0.066 | 0.119 | 0.214 | 0.374 |
| $(d_1, d_2) = (3,3)$ | sign | 0.043 | 0.054 | 0.087 | 0.183 | 0.294 | 0.049 | 0.056 | 0.096 | 0.169 | 0.291 | 0.045 | 0.062 | 0.093 | 0.179 | 0.314 |
| | Spearman | 0.042 | 0.069 | 0.099 | 0.217 | 0.369 | 0.049 | 0.075 | 0.129 | 0.212 | 0.382 | 0.047 | 0.065 | 0.125 | 0.239 | 0.421 |
| | Kendall | 0.049 | 0.066 | 0.098 | 0.220 | 0.357 | 0.055 | 0.071 | 0.113 | 0.207 | 0.356 | 0.044 | 0.073 | 0.110 | 0.218 | 0.397 |
| | vdW | 0.045 | 0.066 | 0.105 | 0.219 | 0.373 | 0.048 | 0.077 | 0.133 | 0.214 | 0.383 | 0.039 | 0.064 | 0.126 | 0.239 | 0.439 |
| | dCov | 0.043 | 0.056 | 0.097 | 0.201 | 0.359 | 0.053 | 0.073 | 0.108 | 0.199 | 0.353 | 0.041 | 0.062 | 0.105 | 0.213 | 0.380 |
| | dCvdW | 0.040 | 0.057 | 0.093 | 0.204 | 0.363 | 0.050 | 0.067 | 0.111 | 0.200 | 0.365 | 0.039 | 0.058 | 0.112 | 0.219 | 0.391 |
| | Wilks | 0.052 | 0.064 | 0.122 | 0.277 | 0.462 | 0.058 | 0.084 | 0.145 | 0.263 | 0.448 | 0.045 | 0.057 | 0.135 | 0.259 | 0.463 |
| | SzR | 0.042 | 0.060 | 0.119 | 0.249 | 0.386 | 0.045 | 0.066 | 0.123 | 0.199 | 0.379 | 0.047 | 0.070 | 0.099 | 0.195 | 0.375 |
| $(d_1, d_2) = (5,5)$ | sign | 0.047 | 0.060 | 0.082 | 0.143 | 0.244 | 0.049 | 0.063 | 0.081 | 0.151 | 0.259 | 0.047 | 0.058 | 0.104 | 0.150 | 0.328 |
| | Spearman | 0.055 | 0.067 | 0.086 | 0.159 | 0.283 | 0.046 | 0.067 | 0.093 | 0.167 | 0.291 | 0.050 | 0.059 | 0.096 | 0.166 | 0.341 |
| | Kendall | 0.054 | 0.069 | 0.087 | 0.159 | 0.279 | 0.054 | 0.067 | 0.094 | 0.173 | 0.279 | 0.047 | 0.053 | 0.113 | 0.172 | 0.321 |
| | vdW | 0.049 | 0.065 | 0.091 | 0.163 | 0.308 | 0.051 | 0.071 | 0.097 | 0.185 | 0.315 | 0.043 | 0.061 | 0.113 | 0.180 | 0.377 |
| | dCov | 0.053 | 0.061 | 0.087 | 0.153 | 0.258 | 0.045 | 0.061 | 0.097 | 0.167 | 0.281 | 0.048 | 0.063 | 0.099 | 0.165 | 0.314 |
| | dCvdW | 0.055 | 0.060 | 0.095 | 0.169 | 0.295 | 0.049 | 0.063 | 0.092 | 0.177 | 0.299 | 0.048 | 0.059 | 0.115 | 0.177 | 0.365 |
| | Wilks | 0.059 | 0.066 | 0.133 | 0.291 | 0.528 | 0.056 | 0.069 | 0.133 | 0.272 | 0.491 | 0.041 | 0.061 | 0.121 | 0.257 | 0.502 |
| | SzR | 0.047 | 0.070 | 0.113 | 0.253 | 0.410 | 0.037 | 0.050 | 0.109 | 0.215 | 0.406 | 0.038 | 0.067 | 0.126 | 0.221 | 0.397 |
| $(d_1, d_2) = (7,7)$ | sign | 0.047 | 0.057 | 0.085 | 0.097 | 0.194 | 0.051 | 0.053 | 0.097 | 0.127 | 0.237 | 0.045 | 0.066 | 0.073 | 0.149 | 0.257 |
| | Spearman | 0.055 | 0.076 | 0.091 | 0.109 | 0.209 | 0.057 | 0.070 | 0.085 | 0.129 | 0.227 | 0.061 | 0.062 | 0.086 | 0.142 | 0.267 |
| | Kendall | 0.055 | 0.064 | 0.089 | 0.113 | 0.203 | 0.053 | 0.064 | 0.087 | 0.131 | 0.226 | 0.051 | 0.059 | 0.083 | 0.149 | 0.260 |
| | vdW | 0.051 | 0.062 | 0.085 | 0.113 | 0.225 | 0.054 | 0.057 | 0.095 | 0.139 | 0.265 | 0.051 | 0.059 | 0.078 | 0.165 | 0.313 |
| | dCov | 0.045 | 0.065 | 0.071 | 0.092 | 0.173 | 0.044 | 0.053 | 0.079 | 0.125 | 0.191 | 0.043 | 0.052 | 0.063 | 0.121 | 0.248 |
| | dCvdW | 0.050 | 0.062 | 0.089 | 0.108 | 0.215 | 0.051 | 0.058 | 0.093 | 0.134 | 0.263 | 0.049 | 0.066 | 0.077 | 0.163 | 0.299 |
| | Wilks | 0.065 | 0.081 | 0.167 | 0.299 | 0.553 | 0.053 | 0.068 | 0.131 | 0.281 | 0.559 | 0.051 | 0.073 | 0.130 | 0.270 | 0.538 |
| | SzR | 0.045 | 0.072 | 0.130 | 0.213 | 0.448 | 0.048 | 0.065 | 0.115 | 0.218 | 0.421 | 0.048 | 0.061 | 0.122 | 0.219 | 0.411 |

Table S.4: Empirical powers, under the heavy-tailed elliptical $t_2$ (case (f)), of *(i)* the center-outward sign (quadrant) test, *(ii)* the center-outward Spearman test, *(iii)* the center-outward Kendall test, *(iv)* the center-outward van der Waerden test, *(v)* the center-outward Wilcoxon score version of the distance covariance test ("dCov"), *(vi)* the center-outward van der Waerden score version of the distance covariance test ("dCvdW"), *(vii)* Wilks' likelihood ratio test, and *(viii)* Székely and Rizzo's permutation-based distance covariance test. Based on 1,000 replications of a sample with dimensions $d_1 = d_2 \in \{2, 3, 5, 7\}$, size $n \in \{432, 864, 1728\}$, and parameter values $\delta^{(n)} = n^{-1/2}\tau$, $\tau \in \{0, 0.2, 0.4, 0.6, 0.8\}$.

| | $n$ | 432 | | | | | 864 | | | | | 1728 | | | | |
|---|---|---|---|---|---|---|---|---|---|---|---|---|---|---|---|---|
| | $\tau$ | 0 | 0.2 | 0.4 | 0.6 | 0.8 | 0 | 0.2 | 0.4 | 0.6 | 0.8 | 0 | 0.2 | 0.4 | 0.6 | 0.8 |
| $(d_1, d_2) = (2, 2)$ | sign | 0.046 | 0.084 | 0.147 | 0.334 | 0.555 | 0.047 | 0.083 | 0.168 | 0.347 | 0.608 | 0.050 | 0.067 | 0.167 | 0.352 | 0.586 |
| | Spearman | 0.047 | 0.084 | 0.205 | 0.445 | 0.695 | 0.047 | 0.075 | 0.206 | 0.427 | 0.711 | 0.047 | 0.077 | 0.201 | 0.437 | 0.716 |
| | Kendall | 0.047 | 0.089 | 0.205 | 0.445 | 0.697 | 0.047 | 0.083 | 0.211 | 0.433 | 0.709 | 0.053 | 0.083 | 0.215 | 0.429 | 0.725 |
| | vdW | 0.048 | 0.087 | 0.208 | 0.469 | 0.738 | 0.049 | 0.081 | 0.216 | 0.471 | 0.740 | 0.049 | 0.081 | 0.221 | 0.481 | 0.761 |
| | dCov | 0.043 | 0.083 | 0.201 | 0.424 | 0.676 | 0.051 | 0.078 | 0.187 | 0.421 | 0.688 | 0.045 | 0.069 | 0.191 | 0.413 | 0.707 |
| | dCvdW | 0.047 | 0.085 | 0.203 | 0.438 | 0.697 | 0.045 | 0.085 | 0.203 | 0.443 | 0.709 | 0.046 | 0.079 | 0.206 | 0.437 | 0.739 |
| | Wilks | 0.070 | 0.089 | 0.143 | 0.278 | 0.509 | 0.063 | 0.083 | 0.139 | 0.262 | 0.482 | 0.061 | 0.079 | 0.126 | 0.248 | 0.467 |
| | SzR | 0.044 | 0.062 | 0.163 | 0.359 | 0.667 | 0.048 | 0.077 | 0.165 | 0.381 | 0.679 | 0.049 | 0.077 | 0.190 | 0.394 | 0.700 |
| $(d_1, d_2) = (3, 3)$ | sign | 0.041 | 0.059 | 0.173 | 0.335 | 0.635 | 0.052 | 0.080 | 0.175 | 0.404 | 0.667 | 0.046 | 0.077 | 0.185 | 0.393 | 0.675 |
| | Spearman | 0.052 | 0.082 | 0.185 | 0.387 | 0.708 | 0.051 | 0.073 | 0.192 | 0.409 | 0.673 | 0.052 | 0.076 | 0.200 | 0.421 | 0.701 |
| | Kendall | 0.059 | 0.077 | 0.199 | 0.401 | 0.701 | 0.048 | 0.073 | 0.193 | 0.440 | 0.719 | 0.047 | 0.073 | 0.197 | 0.433 | 0.719 |
| | vdW | 0.056 | 0.075 | 0.203 | 0.439 | 0.747 | 0.054 | 0.076 | 0.215 | 0.489 | 0.757 | 0.044 | 0.077 | 0.233 | 0.490 | 0.784 |
| | dCov | 0.059 | 0.077 | 0.175 | 0.373 | 0.683 | 0.047 | 0.077 | 0.187 | 0.407 | 0.673 | 0.052 | 0.073 | 0.191 | 0.420 | 0.683 |
| | dCvdW | 0.059 | 0.069 | 0.201 | 0.414 | 0.723 | 0.044 | 0.079 | 0.214 | 0.458 | 0.744 | 0.051 | 0.073 | 0.214 | 0.463 | 0.754 |
| | Wilks | 0.091 | 0.117 | 0.164 | 0.306 | 0.533 | 0.073 | 0.103 | 0.152 | 0.255 | 0.501 | 0.070 | 0.091 | 0.142 | 0.263 | 0.465 |
| | SzR | 0.042 | 0.060 | 0.132 | 0.285 | 0.637 | 0.046 | 0.057 | 0.159 | 0.352 | 0.655 | 0.046 | 0.067 | 0.169 | 0.348 | 0.726 |
| $(d_1, d_2) = (5, 5)$ | sign | 0.049 | 0.075 | 0.128 | 0.318 | 0.577 | 0.037 | 0.077 | 0.160 | 0.345 | 0.651 | 0.053 | 0.076 | 0.167 | 0.402 | 0.709 |
| | Spearman | 0.057 | 0.084 | 0.135 | 0.279 | 0.548 | 0.050 | 0.069 | 0.152 | 0.312 | 0.583 | 0.046 | 0.060 | 0.161 | 0.336 | 0.624 |
| | Kendall | 0.051 | 0.079 | 0.136 | 0.297 | 0.570 | 0.057 | 0.069 | 0.162 | 0.328 | 0.628 | 0.055 | 0.073 | 0.161 | 0.363 | 0.665 |
| | vdW | 0.047 | 0.083 | 0.146 | 0.339 | 0.649 | 0.046 | 0.077 | 0.186 | 0.400 | 0.727 | 0.053 | 0.076 | 0.190 | 0.454 | 0.779 |
| | dCov | 0.060 | 0.078 | 0.138 | 0.274 | 0.541 | 0.048 | 0.070 | 0.157 | 0.307 | 0.583 | 0.051 | 0.076 | 0.154 | 0.328 | 0.631 |
| | dCvdW | 0.053 | 0.081 | 0.154 | 0.335 | 0.638 | 0.045 | 0.077 | 0.187 | 0.375 | 0.709 | 0.055 | 0.071 | 0.186 | 0.436 | 0.757 |
| | Wilks | 0.127 | 0.160 | 0.223 | 0.326 | 0.601 | 0.127 | 0.139 | 0.206 | 0.313 | 0.529 | 0.112 | 0.125 | 0.172 | 0.291 | 0.479 |
| | SzR | 0.051 | 0.056 | 0.100 | 0.234 | 0.499 | 0.053 | 0.057 | 0.121 | 0.270 | 0.543 | 0.050 | 0.066 | 0.121 | 0.281 | 0.637 |
| $(d_1, d_2) = (7, 7)$ | sign | 0.050 | 0.059 | 0.105 | 0.208 | 0.429 | 0.049 | 0.061 | 0.119 | 0.261 | 0.545 | 0.050 | 0.079 | 0.134 | 0.323 | 0.616 |
| | Spearman | 0.062 | 0.073 | 0.105 | 0.208 | 0.365 | 0.050 | 0.075 | 0.121 | 0.242 | 0.448 | 0.052 | 0.072 | 0.120 | 0.261 | 0.503 |
| | Kendall | 0.061 | 0.071 | 0.106 | 0.208 | 0.398 | 0.054 | 0.065 | 0.125 | 0.261 | 0.491 | 0.049 | 0.076 | 0.128 | 0.271 | 0.516 |
| | vdW | 0.058 | 0.063 | 0.111 | 0.238 | 0.463 | 0.047 | 0.073 | 0.135 | 0.313 | 0.592 | 0.045 | 0.076 | 0.163 | 0.359 | 0.668 |
| | dCov | 0.047 | 0.054 | 0.090 | 0.175 | 0.323 | 0.042 | 0.065 | 0.102 | 0.205 | 0.405 | 0.049 | 0.056 | 0.107 | 0.227 | 0.455 |
| | dCvdW | 0.051 | 0.059 | 0.115 | 0.229 | 0.451 | 0.045 | 0.064 | 0.126 | 0.305 | 0.572 | 0.046 | 0.083 | 0.155 | 0.358 | 0.643 |
| | Wilks | 0.142 | 0.188 | 0.231 | 0.381 | 0.614 | 0.139 | 0.158 | 0.198 | 0.331 | 0.538 | 0.133 | 0.168 | 0.203 | 0.296 | 0.479 |
| | SzR | 0.044 | 0.056 | 0.092 | 0.188 | 0.373 | 0.045 | 0.060 | 0.119 | 0.207 | 0.438 | 0.046 | 0.073 | 0.110 | 0.252 | 0.499 |

Table S.5: Empirical powers, under the heavy-tailed elliptical Cauchy distribution (g), of *(i)* the center-outward sign (quadrant) test, *(ii)* the center-outward Spearman test, *(iii)* the center-outward Kendall test, *(iv)* the center-outward van der Waerden test, *(v)* the center-outward Wilcoxon score version of the distance covariance test ("dCov"), *(vi)* the center-outward van der Waerden score version of the distance covariance test ("dCvdW"), *(vii)* Wilks' likelihood ratio test, and *(viii)* Székely and Rizzo's permutation-based distance covariance test. Based on 1,000 replications of a sample with dimensions $d_1 = d_2 \in \{2, 3, 5, 7\}$, size $n \in \{432, 864, 1728\}$, and parameter values $\delta^{(n)} = n^{-1/2}\tau$, $\tau \in \{0, 0.2, 0.4, 0.6, 0.8\}$.

| | $n$ | 432 | | | | | 864 | | | | | 1728 | | | | |
|---|---|---|---|---|---|---|---|---|---|---|---|---|---|---|---|---|
| | $\tau$ | 0 | 0.2 | 0.4 | 0.6 | 0.8 | 0 | 0.2 | 0.4 | 0.6 | 0.8 | 0 | 0.2 | 0.4 | 0.6 | 0.8 |
| $(d_1, d_2) = (2, 2)$ | sign | 0.047 | 0.179 | 0.517 | 0.819 | 0.942 | 0.042 | 0.204 | 0.578 | 0.864 | 0.965 | 0.052 | 0.218 | 0.613 | 0.905 | 0.967 |
| | Spearman | 0.040 | 0.293 | 0.766 | 0.966 | 0.992 | 0.049 | 0.311 | 0.833 | 0.977 | 0.993 | 0.045 | 0.358 | 0.863 | 0.984 | 0.986 |
| | Kendall | 0.049 | 0.277 | 0.740 | 0.948 | 0.990 | 0.051 | 0.285 | 0.794 | 0.969 | 0.992 | 0.049 | 0.331 | 0.817 | 0.979 | 0.984 |
| | vdW | 0.041 | 0.367 | 0.833 | 0.981 | 0.994 | 0.049 | 0.397 | 0.897 | 0.988 | 0.994 | 0.043 | 0.466 | 0.930 | 0.987 | 0.991 |
| | dCov | 0.050 | 0.242 | 0.713 | 0.941 | 0.988 | 0.041 | 0.273 | 0.783 | 0.971 | 0.991 | 0.049 | 0.313 | 0.808 | 0.977 | 0.985 |
| | dCvdW | 0.049 | 0.251 | 0.737 | 0.953 | 0.989 | 0.043 | 0.299 | 0.813 | 0.977 | 0.991 | 0.050 | 0.341 | 0.839 | 0.983 | 0.985 |
| | Wilks | 0.053 | 0.242 | 0.424 | 0.594 | 0.817 | 0.045 | 0.215 | 0.397 | 0.581 | 0.805 | 0.041 | 0.221 | 0.397 | 0.577 | 0.767 |
| | SzR | 0.041 | 0.360 | 0.588 | 0.815 | 0.973 | 0.047 | 0.426 | 0.683 | 0.910 | 0.996 | 0.055 | 0.506 | 0.794 | 0.981 | 1.000 |
| $(d_1, d_2) = (3, 3)$ | sign | 0.051 | 0.183 | 0.549 | 0.875 | 0.976 | 0.061 | 0.217 | 0.636 | 0.927 | 0.981 | 0.042 | 0.224 | 0.683 | 0.939 | 0.977 |
| | Spearman | 0.053 | 0.265 | 0.769 | 0.973 | 0.995 | 0.056 | 0.321 | 0.859 | 0.986 | 0.992 | 0.043 | 0.382 | 0.902 | 0.973 | 0.989 |
| | Kendall | 0.050 | 0.244 | 0.714 | 0.967 | 0.994 | 0.061 | 0.277 | 0.807 | 0.977 | 0.990 | 0.047 | 0.333 | 0.848 | 0.969 | 0.987 |
| | vdW | 0.049 | 0.304 | 0.822 | 0.983 | 0.995 | 0.053 | 0.387 | 0.897 | 0.989 | 0.992 | 0.041 | 0.463 | 0.941 | 0.975 | 0.989 |
| | dCov | 0.055 | 0.229 | 0.709 | 0.956 | 0.993 | 0.055 | 0.278 | 0.801 | 0.977 | 0.991 | 0.043 | 0.321 | 0.847 | 0.969 | 0.986 |
| | dCvdW | 0.050 | 0.241 | 0.739 | 0.960 | 0.992 | 0.058 | 0.304 | 0.809 | 0.976 | 0.990 | 0.044 | 0.345 | 0.865 | 0.970 | 0.985 |
| | Wilks | 0.103 | 0.283 | 0.485 | 0.664 | 0.868 | 0.089 | 0.253 | 0.473 | 0.661 | 0.837 | 0.061 | 0.245 | 0.487 | 0.643 | 0.830 |
| | SzR | 0.052 | 0.352 | 0.562 | 0.799 | 0.943 | 0.046 | 0.402 | 0.639 | 0.865 | 0.982 | 0.052 | 0.488 | 0.765 | 0.956 | 0.998 |
| $(d_1, d_2) = (5, 5)$ | sign | 0.042 | 0.131 | 0.431 | 0.824 | 0.963 | 0.049 | 0.174 | 0.554 | 0.881 | 0.975 | 0.042 | 0.187 | 0.625 | 0.922 | 0.972 |
| | Spearman | 0.049 | 0.197 | 0.641 | 0.939 | 0.993 | 0.048 | 0.252 | 0.775 | 0.981 | 0.989 | 0.045 | 0.296 | 0.844 | 0.974 | 0.988 |
| | Kendall | 0.051 | 0.177 | 0.560 | 0.916 | 0.985 | 0.053 | 0.204 | 0.695 | 0.953 | 0.979 | 0.051 | 0.259 | 0.761 | 0.963 | 0.982 |
| | vdW | 0.049 | 0.198 | 0.661 | 0.955 | 0.991 | 0.047 | 0.276 | 0.816 | 0.983 | 0.988 | 0.042 | 0.343 | 0.883 | 0.970 | 0.986 |
| | dCov | 0.053 | 0.170 | 0.566 | 0.916 | 0.985 | 0.048 | 0.202 | 0.704 | 0.968 | 0.984 | 0.041 | 0.256 | 0.772 | 0.967 | 0.985 |
| | dCvdW | 0.048 | 0.163 | 0.575 | 0.923 | 0.985 | 0.050 | 0.217 | 0.715 | 0.962 | 0.982 | 0.040 | 0.269 | 0.792 | 0.965 | 0.983 |
| | Wilks | 0.155 | 0.369 | 0.577 | 0.757 | 0.919 | 0.133 | 0.340 | 0.541 | 0.736 | 0.891 | 0.095 | 0.330 | 0.509 | 0.705 | 0.873 |
| | SzR | 0.054 | 0.331 | 0.532 | 0.754 | 0.896 | 0.054 | 0.385 | 0.662 | 0.835 | 0.951 | 0.041 | 0.499 | 0.740 | 0.936 | 0.989 |
| $(d_1, d_2) = (7, 7)$ | sign | 0.046 | 0.104 | 0.311 | 0.641 | 0.898 | 0.044 | 0.136 | 0.405 | 0.797 | 0.955 | 0.059 | 0.165 | 0.539 | 0.878 | 0.971 |
| | Spearman | 0.058 | 0.156 | 0.489 | 0.844 | 0.976 | 0.055 | 0.196 | 0.606 | 0.935 | 0.988 | 0.047 | 0.230 | 0.765 | 0.958 | 0.983 |
| | Kendall | 0.055 | 0.134 | 0.411 | 0.777 | 0.956 | 0.053 | 0.163 | 0.511 | 0.904 | 0.983 | 0.061 | 0.191 | 0.651 | 0.935 | 0.978 |
| | vdW | 0.047 | 0.144 | 0.489 | 0.859 | 0.975 | 0.048 | 0.195 | 0.621 | 0.941 | 0.987 | 0.059 | 0.243 | 0.778 | 0.960 | 0.981 |
| | dCov | 0.052 | 0.111 | 0.372 | 0.759 | 0.961 | 0.046 | 0.146 | 0.515 | 0.889 | 0.979 | 0.053 | 0.173 | 0.646 | 0.939 | 0.979 |
| | dCvdW | 0.050 | 0.120 | 0.396 | 0.777 | 0.957 | 0.046 | 0.165 | 0.533 | 0.902 | 0.978 | 0.062 | 0.216 | 0.677 | 0.940 | 0.978 |
| | Wilks | 0.218 | 0.448 | 0.654 | 0.818 | 0.931 | 0.148 | 0.401 | 0.623 | 0.777 | 0.914 | 0.141 | 0.361 | 0.553 | 0.761 | 0.880 |
| | SzR | 0.045 | 0.329 | 0.582 | 0.721 | 0.855 | 0.048 | 0.424 | 0.668 | 0.810 | 0.949 | 0.052 | 0.475 | 0.731 | 0.896 | 0.983 |

Table S.6: Empirical powers, under the heavy-tailed marginal $t_2$ distribution (h), of *(i)* the center-outward sign (quadrant) test, *(ii)* the center-outward Spearman test, *(iii)* the center-outward Kendall test, *(iv)* the center-outward van der Waerden test, *(v)* the center-outward Wilcoxon score version of the distance covariance test ("dCov"), *(vi)* the center-outward van der Waerden score version of the distance covariance test ("dCvdW"), *(vii)* Wilks' likelihood ratio test, and *(viii)* Székely and Rizzo's permutation-based distance covariance test. Based on 1,000 replications of a sample with dimensions $d_1 = d_2 \in \{2, 3, 5, 7\}$, size $n \in \{432, 864, 1728\}$, and parameter values $\delta^{(n)} = n^{-1/2}\tau$, $\tau \in \{0, 0.2, 0.4, 0.6, 0.8\}$.

|  |  | $n$ | 432 | | | | | 864 | | | | | 1728 | | | | |
|---|---|---|---|---|---|---|---|---|---|---|---|---|---|---|---|---|---|
|  |  | $\tau$ | 0 | 0.2 | 0.4 | 0.6 | 0.8 | 0 | 0.2 | 0.4 | 0.6 | 0.8 | 0 | 0.2 | 0.4 | 0.6 | 0.8 |
| $(d_1,d_2)=(2,2)$ | | sign | 0.045 | 0.081 | 0.145 | 0.309 | 0.515 | 0.055 | 0.066 | 0.151 | 0.308 | 0.537 | 0.055 | 0.065 | 0.149 | 0.332 | 0.535 |
| | | Spearman | 0.049 | 0.082 | 0.185 | 0.397 | 0.679 | 0.053 | 0.074 | 0.205 | 0.423 | 0.702 | 0.040 | 0.075 | 0.203 | 0.432 | 0.689 |
| | | Kendall | 0.041 | 0.082 | 0.195 | 0.421 | 0.676 | 0.053 | 0.071 | 0.202 | 0.425 | 0.686 | 0.041 | 0.077 | 0.204 | 0.428 | 0.691 |
| | | vdW | 0.043 | 0.087 | 0.199 | 0.434 | 0.702 | 0.054 | 0.087 | 0.217 | 0.452 | 0.732 | 0.048 | 0.083 | 0.211 | 0.459 | 0.739 |
| | | dCov | 0.048 | 0.075 | 0.165 | 0.385 | 0.648 | 0.055 | 0.080 | 0.183 | 0.404 | 0.661 | 0.049 | 0.070 | 0.205 | 0.411 | 0.663 |
| | | dCvdW | 0.041 | 0.084 | 0.175 | 0.401 | 0.667 | 0.054 | 0.079 | 0.192 | 0.427 | 0.677 | 0.047 | 0.072 | 0.201 | 0.432 | 0.682 |
| | | Wilks | 0.062 | 0.087 | 0.140 | 0.295 | 0.491 | 0.057 | 0.082 | 0.139 | 0.267 | 0.452 | 0.055 | 0.067 | 0.126 | 0.231 | 0.431 |
| | | SzR | 0.048 | 0.058 | 0.167 | 0.379 | 0.639 | 0.051 | 0.076 | 0.159 | 0.365 | 0.645 | 0.041 | 0.077 | 0.175 | 0.432 | 0.675 |
| $(d_1,d_2)=(3,3)$ | | sign | 0.053 | 0.067 | 0.151 | 0.337 | 0.585 | 0.055 | 0.079 | 0.159 | 0.335 | 0.628 | 0.051 | 0.079 | 0.163 | 0.351 | 0.634 |
| | | Spearman | 0.056 | 0.082 | 0.183 | 0.378 | 0.652 | 0.048 | 0.081 | 0.191 | 0.390 | 0.703 | 0.049 | 0.079 | 0.195 | 0.439 | 0.724 |
| | | Kendall | 0.055 | 0.078 | 0.180 | 0.359 | 0.655 | 0.051 | 0.079 | 0.183 | 0.403 | 0.711 | 0.044 | 0.079 | 0.191 | 0.415 | 0.719 |
| | | vdW | 0.055 | 0.080 | 0.201 | 0.406 | 0.706 | 0.046 | 0.076 | 0.200 | 0.435 | 0.766 | 0.047 | 0.085 | 0.214 | 0.476 | 0.789 |
| | | dCov | 0.047 | 0.078 | 0.173 | 0.351 | 0.625 | 0.056 | 0.080 | 0.177 | 0.367 | 0.687 | 0.041 | 0.078 | 0.186 | 0.396 | 0.701 |
| | | dCvdW | 0.051 | 0.076 | 0.189 | 0.380 | 0.677 | 0.051 | 0.082 | 0.187 | 0.403 | 0.741 | 0.039 | 0.089 | 0.194 | 0.437 | 0.746 |
| | | Wilks | 0.089 | 0.105 | 0.150 | 0.296 | 0.543 | 0.087 | 0.093 | 0.140 | 0.281 | 0.512 | 0.081 | 0.099 | 0.139 | 0.260 | 0.473 |
| | | SzR | 0.053 | 0.061 | 0.138 | 0.309 | 0.631 | 0.039 | 0.064 | 0.170 | 0.353 | 0.655 | 0.046 | 0.073 | 0.162 | 0.360 | 0.696 |
| $(d_1,d_2)=(5,5)$ | | sign | 0.043 | 0.066 | 0.121 | 0.277 | 0.502 | 0.053 | 0.071 | 0.128 | 0.283 | 0.569 | 0.051 | 0.065 | 0.159 | 0.347 | 0.625 |
| | | Spearman | 0.048 | 0.071 | 0.127 | 0.277 | 0.511 | 0.049 | 0.070 | 0.147 | 0.287 | 0.561 | 0.053 | 0.084 | 0.166 | 0.349 | 0.639 |
| | | Kendall | 0.052 | 0.067 | 0.119 | 0.266 | 0.471 | 0.055 | 0.071 | 0.131 | 0.287 | 0.533 | 0.048 | 0.073 | 0.154 | 0.328 | 0.599 |
| | | vdW | 0.055 | 0.070 | 0.145 | 0.319 | 0.594 | 0.048 | 0.079 | 0.158 | 0.348 | 0.645 | 0.055 | 0.091 | 0.189 | 0.412 | 0.735 |
| | | dCov | 0.045 | 0.067 | 0.115 | 0.265 | 0.487 | 0.055 | 0.075 | 0.143 | 0.276 | 0.548 | 0.055 | 0.087 | 0.155 | 0.332 | 0.609 |
| | | dCvdW | 0.047 | 0.070 | 0.133 | 0.303 | 0.573 | 0.052 | 0.074 | 0.148 | 0.328 | 0.624 | 0.049 | 0.076 | 0.183 | 0.397 | 0.708 |
| | | Wilks | 0.113 | 0.127 | 0.199 | 0.341 | 0.580 | 0.107 | 0.135 | 0.175 | 0.279 | 0.534 | 0.107 | 0.109 | 0.161 | 0.281 | 0.501 |
| | | SzR | 0.047 | 0.054 | 0.116 | 0.258 | 0.537 | 0.048 | 0.060 | 0.128 | 0.283 | 0.590 | 0.054 | 0.073 | 0.124 | 0.317 | 0.664 |
| $(d_1,d_2)=(7,7)$ | | sign | 0.051 | 0.061 | 0.094 | 0.179 | 0.369 | 0.042 | 0.056 | 0.115 | 0.226 | 0.448 | 0.049 | 0.059 | 0.133 | 0.257 | 0.529 |
| | | Spearman | 0.065 | 0.073 | 0.107 | 0.163 | 0.332 | 0.057 | 0.069 | 0.105 | 0.229 | 0.406 | 0.057 | 0.063 | 0.105 | 0.244 | 0.479 |
| | | Kendall | 0.058 | 0.066 | 0.098 | 0.160 | 0.303 | 0.043 | 0.065 | 0.097 | 0.204 | 0.377 | 0.055 | 0.062 | 0.117 | 0.216 | 0.427 |
| | | vdW | 0.051 | 0.065 | 0.107 | 0.197 | 0.406 | 0.039 | 0.067 | 0.120 | 0.262 | 0.505 | 0.059 | 0.065 | 0.137 | 0.302 | 0.593 |
| | | dCov | 0.046 | 0.063 | 0.090 | 0.141 | 0.279 | 0.050 | 0.065 | 0.100 | 0.185 | 0.372 | 0.046 | 0.063 | 0.095 | 0.223 | 0.428 |
| | | dCvdW | 0.050 | 0.062 | 0.098 | 0.183 | 0.386 | 0.037 | 0.066 | 0.119 | 0.250 | 0.501 | 0.049 | 0.065 | 0.133 | 0.288 | 0.572 |
| | | Wilks | 0.139 | 0.161 | 0.226 | 0.381 | 0.593 | 0.125 | 0.149 | 0.195 | 0.327 | 0.539 | 0.127 | 0.131 | 0.195 | 0.298 | 0.529 |
| | | SzR | 0.038 | 0.056 | 0.119 | 0.210 | 0.471 | 0.043 | 0.060 | 0.113 | 0.237 | 0.549 | 0.058 | 0.071 | 0.123 | 0.268 | 0.556 |

Table S.7: Empirical powers, under the heavy-tailed marginal Cauchy distribution (i), of *(i)* the center-outward sign (quadrant) test, *(ii)* the center-outward Spearman test, *(iii)* the center-outward Kendall test, *(iv)* the center-outward van der Waerden test, *(v)* the center-outward Wilcoxon score version of the distance covariance test ("dCov"), *(vi)* the center-outward van der Waerden score version of the distance covariance test ("dCvdW"), *(vii)* Wilks' likelihood ratio test, and *(viii)* Székely and Rizzo's permutation-based distance covariance test. Based on 1,000 replications of a sample with dimensions $d_1 = d_2 \in \{2, 3, 5, 7\}$, size $n \in \{432, 864, 1728\}$, and parameter values $\delta^{(n)} = n^{-1/2}\tau$, $\tau \in \{0, 0.2, 0.4, 0.6, 0.8\}$.

|  |  | $n$ | 432 | | | | | 864 | | | | | 1728 | | | | |
|---|---|---|---|---|---|---|---|---|---|---|---|---|---|---|---|---|---|
|  |  | $\tau$ | 0 | 0.2 | 0.4 | 0.6 | 0.8 | 0 | 0.2 | 0.4 | 0.6 | 0.8 | 0 | 0.2 | 0.4 | 0.6 | 0.8 |
| $(d_1, d_2) = (2, 2)$ | | sign | 0.053 | 0.187 | 0.548 | 0.857 | 0.973 | 0.043 | 0.229 | 0.618 | 0.881 | 0.974 | 0.047 | 0.241 | 0.665 | 0.905 | 0.968 |
| | | Spearman | 0.050 | 0.344 | 0.828 | 0.975 | 0.993 | 0.042 | 0.375 | 0.880 | 0.983 | 0.992 | 0.047 | 0.426 | 0.905 | 0.979 | 0.985 |
| | | Kendall | 0.052 | 0.291 | 0.755 | 0.955 | 0.993 | 0.045 | 0.319 | 0.805 | 0.973 | 0.989 | 0.052 | 0.350 | 0.834 | 0.969 | 0.980 |
| | | vdW | 0.053 | 0.404 | 0.875 | 0.988 | 0.995 | 0.039 | 0.456 | 0.922 | 0.988 | 0.993 | 0.046 | 0.519 | 0.945 | 0.983 | 0.987 |
| | | dCov | 0.055 | 0.295 | 0.755 | 0.967 | 0.992 | 0.037 | 0.315 | 0.825 | 0.974 | 0.991 | 0.051 | 0.355 | 0.854 | 0.973 | 0.982 |
| | | dCvdW | 0.051 | 0.311 | 0.785 | 0.972 | 0.992 | 0.040 | 0.338 | 0.851 | 0.979 | 0.991 | 0.048 | 0.385 | 0.877 | 0.974 | 0.983 |
| | | Wilks | 0.062 | 0.224 | 0.381 | 0.592 | 0.782 | 0.049 | 0.220 | 0.401 | 0.571 | 0.775 | 0.041 | 0.201 | 0.390 | 0.567 | 0.713 |
| | | SzR | 0.040 | 0.334 | 0.580 | 0.831 | 0.971 | 0.059 | 0.413 | 0.715 | 0.923 | 0.995 | 0.046 | 0.510 | 0.817 | 0.975 | 1.000 |
| $(d_1, d_2) = (3, 3)$ | | sign | 0.040 | 0.194 | 0.594 | 0.909 | 0.985 | 0.053 | 0.237 | 0.687 | 0.941 | 0.981 | 0.063 | 0.284 | 0.752 | 0.943 | 0.971 |
| | | Spearman | 0.051 | 0.325 | 0.824 | 0.984 | 0.994 | 0.055 | 0.396 | 0.908 | 0.987 | 0.989 | 0.055 | 0.473 | 0.937 | 0.976 | 0.982 |
| | | Kendall | 0.047 | 0.219 | 0.698 | 0.963 | 0.992 | 0.055 | 0.274 | 0.787 | 0.970 | 0.985 | 0.050 | 0.327 | 0.833 | 0.963 | 0.975 |
| | | vdW | 0.045 | 0.352 | 0.860 | 0.987 | 0.994 | 0.055 | 0.447 | 0.926 | 0.987 | 0.989 | 0.057 | 0.531 | 0.951 | 0.978 | 0.983 |
| | | dCov | 0.058 | 0.253 | 0.777 | 0.980 | 0.994 | 0.051 | 0.329 | 0.858 | 0.983 | 0.986 | 0.050 | 0.405 | 0.899 | 0.971 | 0.979 |
| | | dCvdW | 0.050 | 0.266 | 0.783 | 0.979 | 0.994 | 0.049 | 0.342 | 0.865 | 0.984 | 0.986 | 0.061 | 0.417 | 0.902 | 0.969 | 0.980 |
| | | Wilks | 0.091 | 0.277 | 0.477 | 0.681 | 0.826 | 0.073 | 0.261 | 0.445 | 0.641 | 0.816 | 0.064 | 0.243 | 0.437 | 0.611 | 0.798 |
| | | SzR | 0.065 | 0.352 | 0.576 | 0.795 | 0.939 | 0.061 | 0.416 | 0.682 | 0.903 | 0.983 | 0.039 | 0.477 | 0.791 | 0.960 | 1.000 |
| $(d_1, d_2) = (5, 5)$ | | sign | 0.049 | 0.164 | 0.535 | 0.877 | 0.973 | 0.033 | 0.223 | 0.667 | 0.930 | 0.968 | 0.045 | 0.238 | 0.746 | 0.937 | 0.959 |
| | | Spearman | 0.049 | 0.247 | 0.732 | 0.968 | 0.989 | 0.065 | 0.335 | 0.869 | 0.979 | 0.985 | 0.048 | 0.429 | 0.907 | 0.969 | 0.979 |
| | | Kendall | 0.051 | 0.151 | 0.474 | 0.828 | 0.969 | 0.052 | 0.201 | 0.590 | 0.903 | 0.971 | 0.048 | 0.225 | 0.704 | 0.939 | 0.961 |
| | | vdW | 0.049 | 0.241 | 0.762 | 0.976 | 0.989 | 0.050 | 0.355 | 0.881 | 0.978 | 0.983 | 0.049 | 0.445 | 0.917 | 0.966 | 0.977 |
| | | dCov | 0.050 | 0.193 | 0.656 | 0.945 | 0.989 | 0.055 | 0.272 | 0.796 | 0.974 | 0.978 | 0.058 | 0.363 | 0.875 | 0.967 | 0.974 |
| | | dCvdW | 0.046 | 0.205 | 0.673 | 0.953 | 0.987 | 0.045 | 0.290 | 0.797 | 0.969 | 0.976 | 0.048 | 0.354 | 0.863 | 0.961 | 0.969 |
| | | Wilks | 0.160 | 0.362 | 0.575 | 0.750 | 0.897 | 0.134 | 0.335 | 0.527 | 0.722 | 0.858 | 0.089 | 0.311 | 0.503 | 0.705 | 0.847 |
| | | SzR | 0.039 | 0.341 | 0.567 | 0.739 | 0.871 | 0.045 | 0.389 | 0.631 | 0.825 | 0.958 | 0.048 | 0.488 | 0.745 | 0.936 | 0.991 |
| $(d_1, d_2) = (7, 7)$ | | sign | 0.052 | 0.119 | 0.381 | 0.717 | 0.929 | 0.044 | 0.149 | 0.544 | 0.858 | 0.963 | 0.049 | 0.197 | 0.659 | 0.918 | 0.951 |
| | | Spearman | 0.068 | 0.193 | 0.561 | 0.884 | 0.976 | 0.056 | 0.233 | 0.743 | 0.965 | 0.983 | 0.053 | 0.317 | 0.849 | 0.962 | 0.967 |
| | | Kendall | 0.059 | 0.125 | 0.289 | 0.571 | 0.844 | 0.057 | 0.125 | 0.385 | 0.734 | 0.936 | 0.055 | 0.143 | 0.487 | 0.845 | 0.933 |
| | | vdW | 0.052 | 0.176 | 0.563 | 0.882 | 0.980 | 0.047 | 0.231 | 0.747 | 0.961 | 0.980 | 0.055 | 0.321 | 0.847 | 0.955 | 0.965 |
| | | dCov | 0.044 | 0.129 | 0.428 | 0.803 | 0.971 | 0.042 | 0.166 | 0.635 | 0.929 | 0.977 | 0.057 | 0.233 | 0.756 | 0.953 | 0.962 |
| | | dCvdW | 0.047 | 0.144 | 0.485 | 0.833 | 0.968 | 0.043 | 0.189 | 0.669 | 0.927 | 0.973 | 0.053 | 0.266 | 0.780 | 0.946 | 0.961 |
| | | Wilks | 0.201 | 0.425 | 0.643 | 0.829 | 0.904 | 0.170 | 0.391 | 0.575 | 0.770 | 0.899 | 0.128 | 0.363 | 0.566 | 0.733 | 0.869 |
| | | SzR | 0.053 | 0.356 | 0.528 | 0.750 | 0.884 | 0.048 | 0.417 | 0.649 | 0.838 | 0.950 | 0.043 | 0.473 | 0.714 | 0.915 | 0.988 |

Table S.8: Empirical powers, under the Gaussian mixture alternative (j), of *(i)* the center-outward sign (quadrant) test, *(ii)* the center-outward Spearman test, *(iii)* the center-outward Kendall test, *(iv)* the center-outward van der Waerden test, *(v)* the center-outward Wilcoxon score version of the distance covariance test ("dCov"), *(vi)* the center-outward van der Waerden score version of the distance covariance test ("dCvdW"), *(vii)* Wilks' likelihood ratio test, and *(viii)* Székely and Rizzo's permutation-based distance covariance test. Based on 1,000 replications of a sample with dimensions $d_1 = d_2 \in \{2, 3, 5, 7\}$, size $n \in \{432, 864, 1728\}$, and parameter values $\delta^{(n)} = n^{-1/2}\tau$, $\tau \in \{0, 1, 2, 3, 4\}$.

| | $n$ | 432 | | | | | 864 | | | | | 1728 | | | | |
|---|---|---|---|---|---|---|---|---|---|---|---|---|---|---|---|---|
| | $\tau$ | 0 | 1 | 2 | 3 | 4 | 0 | 1 | 2 | 3 | 4 | 0 | 1 | 2 | 3 | 4 |
| $(d_1, d_2) = (2, 2)$ | sign | 0.062 | 0.063 | 0.123 | 0.245 | 0.409 | 0.049 | 0.065 | 0.137 | 0.234 | 0.424 | 0.052 | 0.062 | 0.118 | 0.263 | 0.401 |
| | Spearman | 0.042 | 0.060 | 0.164 | 0.347 | 0.561 | 0.052 | 0.079 | 0.177 | 0.342 | 0.593 | 0.039 | 0.080 | 0.170 | 0.356 | 0.571 |
| | Kendall | 0.041 | 0.062 | 0.174 | 0.342 | 0.568 | 0.060 | 0.073 | 0.183 | 0.332 | 0.567 | 0.048 | 0.073 | 0.156 | 0.350 | 0.562 |
| | vdW | 0.046 | 0.063 | 0.172 | 0.358 | 0.564 | 0.061 | 0.076 | 0.172 | 0.351 | 0.602 | 0.038 | 0.072 | 0.174 | 0.363 | 0.594 |
| | dCov | 0.045 | 0.057 | 0.159 | 0.322 | 0.532 | 0.051 | 0.072 | 0.167 | 0.323 | 0.555 | 0.044 | 0.068 | 0.173 | 0.345 | 0.547 |
| | dCvdW | 0.041 | 0.055 | 0.171 | 0.330 | 0.529 | 0.051 | 0.073 | 0.166 | 0.325 | 0.564 | 0.042 | 0.074 | 0.174 | 0.345 | 0.554 |
| | Wilks | 0.042 | 0.069 | 0.194 | 0.392 | 0.626 | 0.058 | 0.081 | 0.192 | 0.367 | 0.624 | 0.036 | 0.076 | 0.173 | 0.373 | 0.605 |
| | SzR | 0.048 | 0.061 | 0.174 | 0.352 | 0.580 | 0.042 | 0.068 | 0.168 | 0.330 | 0.582 | 0.039 | 0.072 | 0.175 | 0.347 | 0.567 |
| $(d_1, d_2) = (3, 3)$ | sign | 0.040 | 0.063 | 0.125 | 0.236 | 0.445 | 0.056 | 0.062 | 0.120 | 0.263 | 0.455 | 0.058 | 0.058 | 0.140 | 0.278 | 0.477 |
| | Spearman | 0.051 | 0.063 | 0.128 | 0.309 | 0.549 | 0.057 | 0.073 | 0.154 | 0.336 | 0.557 | 0.051 | 0.065 | 0.168 | 0.357 | 0.575 |
| | Kendall | 0.048 | 0.071 | 0.127 | 0.308 | 0.536 | 0.053 | 0.068 | 0.157 | 0.311 | 0.558 | 0.054 | 0.064 | 0.170 | 0.354 | 0.564 |
| | vdW | 0.042 | 0.066 | 0.136 | 0.313 | 0.552 | 0.056 | 0.068 | 0.155 | 0.334 | 0.578 | 0.061 | 0.064 | 0.174 | 0.374 | 0.606 |
| | dCov | 0.036 | 0.070 | 0.113 | 0.292 | 0.518 | 0.051 | 0.058 | 0.161 | 0.310 | 0.536 | 0.057 | 0.067 | 0.143 | 0.319 | 0.564 |
| | dCvdW | 0.038 | 0.069 | 0.137 | 0.296 | 0.538 | 0.059 | 0.064 | 0.160 | 0.324 | 0.562 | 0.064 | 0.067 | 0.159 | 0.335 | 0.571 |
| | Wilks | 0.041 | 0.081 | 0.174 | 0.393 | 0.657 | 0.060 | 0.071 | 0.191 | 0.382 | 0.650 | 0.053 | 0.066 | 0.177 | 0.408 | 0.648 |
| | SzR | 0.040 | 0.074 | 0.148 | 0.367 | 0.618 | 0.051 | 0.062 | 0.171 | 0.358 | 0.604 | 0.065 | 0.066 | 0.187 | 0.377 | 0.606 |
| $(d_1, d_2) = (5, 5)$ | sign | 0.044 | 0.067 | 0.100 | 0.217 | 0.402 | 0.044 | 0.059 | 0.111 | 0.244 | 0.431 | 0.050 | 0.059 | 0.118 | 0.237 | 0.503 |
| | Spearman | 0.058 | 0.070 | 0.114 | 0.216 | 0.432 | 0.041 | 0.070 | 0.141 | 0.259 | 0.458 | 0.051 | 0.066 | 0.152 | 0.286 | 0.503 |
| | Kendall | 0.059 | 0.061 | 0.108 | 0.242 | 0.455 | 0.043 | 0.052 | 0.129 | 0.247 | 0.454 | 0.055 | 0.063 | 0.137 | 0.282 | 0.513 |
| | vdW | 0.042 | 0.067 | 0.118 | 0.244 | 0.477 | 0.048 | 0.061 | 0.136 | 0.285 | 0.499 | 0.051 | 0.057 | 0.135 | 0.305 | 0.569 |
| | dCov | 0.046 | 0.071 | 0.105 | 0.220 | 0.408 | 0.049 | 0.062 | 0.127 | 0.261 | 0.442 | 0.051 | 0.065 | 0.122 | 0.278 | 0.503 |
| | dCvdW | 0.039 | 0.066 | 0.109 | 0.243 | 0.463 | 0.047 | 0.067 | 0.121 | 0.280 | 0.480 | 0.048 | 0.057 | 0.125 | 0.296 | 0.557 |
| | Wilks | 0.056 | 0.079 | 0.191 | 0.425 | 0.748 | 0.042 | 0.076 | 0.192 | 0.426 | 0.721 | 0.060 | 0.074 | 0.186 | 0.409 | 0.721 |
| | SzR | 0.035 | 0.065 | 0.176 | 0.404 | 0.708 | 0.043 | 0.060 | 0.170 | 0.384 | 0.689 | 0.055 | 0.063 | 0.165 | 0.381 | 0.698 |
| $(d_1, d_2) = (7, 7)$ | sign | 0.038 | 0.056 | 0.086 | 0.169 | 0.305 | 0.049 | 0.057 | 0.110 | 0.210 | 0.355 | 0.049 | 0.076 | 0.113 | 0.220 | 0.424 |
| | Spearman | 0.047 | 0.084 | 0.105 | 0.172 | 0.310 | 0.052 | 0.061 | 0.108 | 0.208 | 0.379 | 0.054 | 0.065 | 0.112 | 0.216 | 0.401 |
| | Kendall | 0.050 | 0.070 | 0.106 | 0.180 | 0.305 | 0.057 | 0.064 | 0.113 | 0.212 | 0.372 | 0.051 | 0.072 | 0.115 | 0.213 | 0.416 |
| | vdW | 0.036 | 0.060 | 0.096 | 0.183 | 0.349 | 0.048 | 0.057 | 0.121 | 0.219 | 0.429 | 0.049 | 0.068 | 0.109 | 0.247 | 0.477 |
| | dCov | 0.035 | 0.054 | 0.079 | 0.146 | 0.263 | 0.047 | 0.059 | 0.096 | 0.177 | 0.326 | 0.044 | 0.047 | 0.094 | 0.169 | 0.354 |
| | dCvdW | 0.028 | 0.060 | 0.105 | 0.182 | 0.340 | 0.051 | 0.056 | 0.116 | 0.215 | 0.399 | 0.054 | 0.066 | 0.114 | 0.219 | 0.458 |
| | Wilks | 0.067 | 0.088 | 0.208 | 0.472 | 0.764 | 0.059 | 0.073 | 0.199 | 0.449 | 0.798 | 0.050 | 0.061 | 0.206 | 0.423 | 0.755 |
| | SzR | 0.049 | 0.068 | 0.172 | 0.427 | 0.723 | 0.057 | 0.063 | 0.180 | 0.414 | 0.769 | 0.045 | 0.064 | 0.178 | 0.404 | 0.744 |

Table S.9: Empirical powers, under the heavy-tailed elliptical mixture alternative (k), of *(i)* the center-outward sign (quadrant) test, *(ii)* the center-outward Spearman test, *(iii)* the center-outward Kendall test, *(iv)* the center-outward van der Waerden test, *(v)* the center-outward Wilcoxon score version of the distance covariance test ("dCov"), *(vi)* the center-outward van der Waerden score version of the distance covariance test ("dCvdW"), *(vii)* Wilks' likelihood ratio test, and *(viii)* Székely and Rizzo's permutation-based distance covariance test. Based on 1,000 replications of a sample with dimensions $d_1 = d_2 \in \{2, 3, 5, 7\}$, size $n \in \{432, 864, 1728\}$, and parameter values $\delta^{(n)} = n^{-1/2}\tau$, $\tau \in \{0, 1, 2, 3, 4\}$.

|  | $n$ | 432 | | | | | 864 | | | | | 1728 | | | | |
|---|---|---|---|---|---|---|---|---|---|---|---|---|---|---|---|---|
|  | $\tau$ | 0 | 1 | 2 | 3 | 4 | 0 | 1 | 2 | 3 | 4 | 0 | 1 | 2 | 3 | 4 |
| $(d_1, d_2) = (2, 2)$ | sign | 0.058 | 0.070 | 0.102 | 0.235 | 0.391 | 0.058 | 0.065 | 0.133 | 0.230 | 0.418 | 0.052 | 0.064 | 0.127 | 0.244 | 0.410 |
|  | Spearman | 0.045 | 0.088 | 0.186 | 0.344 | 0.543 | 0.047 | 0.090 | 0.192 | 0.329 | 0.597 | 0.062 | 0.073 | 0.168 | 0.331 | 0.562 |
|  | Kendall | 0.053 | 0.088 | 0.166 | 0.334 | 0.523 | 0.052 | 0.088 | 0.178 | 0.315 | 0.565 | 0.067 | 0.076 | 0.156 | 0.316 | 0.532 |
|  | vdW | 0.039 | 0.088 | 0.184 | 0.357 | 0.548 | 0.050 | 0.087 | 0.201 | 0.335 | 0.602 | 0.053 | 0.072 | 0.176 | 0.340 | 0.584 |
|  | dCov | 0.047 | 0.079 | 0.158 | 0.331 | 0.517 | 0.047 | 0.091 | 0.173 | 0.310 | 0.572 | 0.063 | 0.067 | 0.165 | 0.316 | 0.556 |
|  | dCvdW | 0.049 | 0.085 | 0.157 | 0.322 | 0.504 | 0.050 | 0.080 | 0.172 | 0.308 | 0.564 | 0.059 | 0.070 | 0.169 | 0.307 | 0.548 |
|  | Wilks | 0.068 | 0.130 | 0.266 | 0.435 | 0.623 | 0.064 | 0.124 | 0.265 | 0.423 | 0.627 | 0.063 | 0.112 | 0.252 | 0.409 | 0.598 |
|  | SzR | 0.045 | 0.116 | 0.264 | 0.491 | 0.717 | 0.050 | 0.098 | 0.248 | 0.449 | 0.713 | 0.060 | 0.081 | 0.210 | 0.430 | 0.683 |
| $(d_1, d_2) = (3, 3)$ | sign | 0.048 | 0.057 | 0.131 | 0.253 | 0.444 | 0.060 | 0.084 | 0.155 | 0.256 | 0.481 | 0.047 | 0.066 | 0.139 | 0.262 | 0.498 |
|  | Spearman | 0.042 | 0.086 | 0.166 | 0.359 | 0.567 | 0.057 | 0.098 | 0.141 | 0.341 | 0.604 | 0.054 | 0.082 | 0.170 | 0.357 | 0.615 |
|  | Kendall | 0.053 | 0.077 | 0.149 | 0.332 | 0.539 | 0.049 | 0.093 | 0.148 | 0.319 | 0.575 | 0.043 | 0.073 | 0.160 | 0.332 | 0.568 |
|  | vdW | 0.041 | 0.077 | 0.158 | 0.352 | 0.561 | 0.059 | 0.089 | 0.159 | 0.345 | 0.600 | 0.055 | 0.073 | 0.155 | 0.356 | 0.624 |
|  | dCov | 0.049 | 0.065 | 0.142 | 0.319 | 0.554 | 0.046 | 0.077 | 0.135 | 0.321 | 0.582 | 0.057 | 0.075 | 0.165 | 0.356 | 0.604 |
|  | dCvdW | 0.048 | 0.071 | 0.145 | 0.322 | 0.524 | 0.051 | 0.079 | 0.152 | 0.323 | 0.570 | 0.047 | 0.072 | 0.155 | 0.345 | 0.594 |
|  | Wilks | 0.061 | 0.156 | 0.301 | 0.511 | 0.696 | 0.068 | 0.179 | 0.286 | 0.500 | 0.686 | 0.084 | 0.123 | 0.265 | 0.447 | 0.680 |
|  | SzR | 0.046 | 0.120 | 0.308 | 0.560 | 0.811 | 0.054 | 0.133 | 0.270 | 0.538 | 0.803 | 0.050 | 0.100 | 0.269 | 0.521 | 0.790 |
| $(d_1, d_2) = (5, 5)$ | sign | 0.043 | 0.069 | 0.099 | 0.199 | 0.382 | 0.040 | 0.062 | 0.111 | 0.249 | 0.413 | 0.050 | 0.055 | 0.109 | 0.237 | 0.470 |
|  | Spearman | 0.059 | 0.076 | 0.150 | 0.262 | 0.482 | 0.056 | 0.080 | 0.154 | 0.314 | 0.512 | 0.051 | 0.075 | 0.133 | 0.312 | 0.565 |
|  | Kendall | 0.055 | 0.069 | 0.145 | 0.240 | 0.455 | 0.056 | 0.068 | 0.126 | 0.273 | 0.448 | 0.058 | 0.058 | 0.120 | 0.258 | 0.527 |
|  | vdW | 0.043 | 0.073 | 0.138 | 0.237 | 0.469 | 0.049 | 0.072 | 0.129 | 0.299 | 0.487 | 0.052 | 0.055 | 0.129 | 0.268 | 0.550 |
|  | dCov | 0.056 | 0.070 | 0.145 | 0.247 | 0.483 | 0.052 | 0.066 | 0.129 | 0.299 | 0.543 | 0.052 | 0.064 | 0.140 | 0.310 | 0.614 |
|  | dCvdW | 0.045 | 0.071 | 0.119 | 0.206 | 0.420 | 0.046 | 0.069 | 0.108 | 0.265 | 0.465 | 0.056 | 0.058 | 0.128 | 0.254 | 0.546 |
|  | Wilks | 0.092 | 0.207 | 0.417 | 0.631 | 0.815 | 0.095 | 0.198 | 0.364 | 0.597 | 0.789 | 0.072 | 0.182 | 0.326 | 0.539 | 0.778 |
|  | SzR | 0.046 | 0.163 | 0.385 | 0.640 | 0.882 | 0.052 | 0.133 | 0.332 | 0.635 | 0.880 | 0.045 | 0.127 | 0.291 | 0.594 | 0.888 |
| $(d_1, d_2) = (7, 7)$ | sign | 0.042 | 0.059 | 0.084 | 0.163 | 0.299 | 0.052 | 0.057 | 0.100 | 0.178 | 0.320 | 0.059 | 0.060 | 0.102 | 0.208 | 0.421 |
|  | Spearman | 0.051 | 0.108 | 0.131 | 0.215 | 0.422 | 0.059 | 0.077 | 0.138 | 0.242 | 0.403 | 0.054 | 0.082 | 0.113 | 0.264 | 0.496 |
|  | Kendall | 0.036 | 0.081 | 0.112 | 0.183 | 0.389 | 0.057 | 0.073 | 0.115 | 0.200 | 0.366 | 0.060 | 0.076 | 0.113 | 0.217 | 0.413 |
|  | vdW | 0.041 | 0.078 | 0.112 | 0.169 | 0.358 | 0.066 | 0.069 | 0.121 | 0.212 | 0.392 | 0.058 | 0.069 | 0.101 | 0.236 | 0.472 |
|  | dCov | 0.034 | 0.075 | 0.100 | 0.179 | 0.429 | 0.053 | 0.057 | 0.119 | 0.244 | 0.458 | 0.064 | 0.070 | 0.124 | 0.289 | 0.582 |
|  | dCvdW | 0.036 | 0.077 | 0.094 | 0.157 | 0.345 | 0.047 | 0.056 | 0.105 | 0.196 | 0.349 | 0.060 | 0.066 | 0.099 | 0.216 | 0.444 |
|  | Wilks | 0.095 | 0.272 | 0.458 | 0.728 | 0.867 | 0.100 | 0.242 | 0.423 | 0.660 | 0.852 | 0.092 | 0.222 | 0.417 | 0.624 | 0.835 |
|  | SzR | 0.045 | 0.170 | 0.385 | 0.699 | 0.904 | 0.039 | 0.149 | 0.380 | 0.709 | 0.914 | 0.060 | 0.150 | 0.364 | 0.682 | 0.920 |



\end{document}